\documentclass{article}

\usepackage{amsmath}

\usepackage[colorlinks,allcolors=blue]{hyperref}

\usepackage[natbib,style=rsl]{biblatex}
\usepackage{csquotes}
\usepackage[french,german,greek,italian,latin,english]{babel}
\usepackage{epigraph}
\usepackage{amssymb,amsfonts,stix}

\def\peirceplus{\mathbin{+\mkern-2mu,}}
\usepackage[T3,T1]{fontenc}
\usepackage[safe]{tipa}
\def\leibnizv{\mbox{\textroundcap{v}}}

\usepackage{pifont,graphicx}
\DeclareUnicodeCharacter{274B}{\pmast}
\newcommand{\fmimp}{\mathrel{\rotatebox[origin=c]{180}{\scriptsize\textbf{C}}}}
\newcommand{\pmthm}{\mathpunct{\text{\scalebox{.5}[1]{$\boldsymbol\vdash$}}}}
\newcommand{\pmast}{\ensuremath{\mathord{\ast}}}
\newcommand{\pmcdot}{\raisebox{.05cm}{$\boldsymbol\cdot$}}
\newcommand{\pmiddf}{\mathrel{=}}
\newcommand{\pmdf}{\quad \text{Df}}
\newcommand{\pmdot}{\mathrel{\hbox{\rule{.3ex}{.3ex}}}}
\newcommand{\pmdott}{\mathrel{\overset{\pmdot}{\pmdot}}}
\newcommand{\pmdottt}{\pmdott\hspace{.1em}\pmdot}

\newcommand{\pmand}{\mathrel{\hbox{\rule{.3ex}{.3ex}}}}

\renewcommand{\lnot}{\mathord{\boldsymbol{\neg}}}
\renewcommand{\land}{\mathbin{\boldsymbol{\wedge}}}
\newcommand{\laddexcl}{\mathbin{\boldsymbol{\barvee}}}

\newcommand{\pmnot}{\mathord{\boldsymbol{\sim}}}
\newcommand{\pmimp}{\mathbin{\boldsymbol{\supset}}}
\newcommand{\pmiff}{\mathbin{\boldsymbol{\equiv}}}
\newcommand{\pmor}{\mathbin{\boldsymbol{\vee}}}
\let\lor\pmor
\newcommand{\pmSome}{\text{\raisebox{.5em}{\rotatebox{180}{\textbf{E}}}}}
\newcommand{\pmbreve}[1]{\mathbf{\breve{\text{$#1$}}}}
\newcommand{\pmcrel}[1]{\pmbreve{#1}} 

\newcommand{\pmcinc}{\mathrel{\boldsymbol{\subset}}}
\newcommand{\pmccap}{\mathrel{\boldsymbol{\cap}}}
\newcommand{\pmccup}{\mathrel{\boldsymbol{\cup}}}
\newcommand{\grassfrown}{\mathrel{\boldsymbol{\smallfrown}}}
\newcommand{\grasssmile}{\mathrel{\boldsymbol{\smallsmile}}}

\newcommand{\pmrcup}{\mathrel{\ooalign{$\hidewidth\raisebox{.1em}{$\boldsymbol{\cdot}$}\hidewidth$\cr$\mathbf{\pmccup}$}}}

\def\peirceor{\mathbin{\ooalign{$\smallsmile$\cr\hfil$\mid$\hfil}}}
\def\ampheck{\mathbin{\rotatebox{90}{$\scurel$}}}
\def\peircearrow{\mathbin{\ooalign{\raisebox{-.2em}{$\boldsymbol\vee$}\cr\hfil$\boldsymbol\mid$\hfil}}}
\def\conseq{\mathbin{\raisebox{\depth}{\rotatebox{180}{\textsc{c}}}}}
\def\subsumed{\mathbin{\ooalign{$\subset$\cr\hidewidth\hbox{$=\mkern.5mu$}\hidewidth}}}

\title{The Genealogy of `$\lor$'}
\author{
    Landon D. C. Elkind\thanks{Much of this research was done while a Killam
    Postdoctoral Fellow at the University of Alberta.} 
    ~(Western Kentucky University)
\and 
    Richard Zach 
    (University of Calgary)
}
\date{}

\begin{document}
\maketitle
\selectlanguage{english}

\begin{abstract}
The use of the symbol~$\lor$ for disjunction in formal logic is
ubiquitous. Where did it come from?  The paper details the evolution
of the symbol~$\lor$ in its historical and logical context. Some
sources say that disjunction in its use as connecting propositions or
formulas was introduced by Peano; others suggest that it originated as
an abbreviation of the Latin word for ``or'', \emph{vel}. We show that
the origin of the symbol~$\lor$ for disjunction can be traced to
Whitehead and Russell's pre-\textit{Principia} work in formal logic.
Because of \textit{Principia}'s influence, its notation was widely
adopted by philosophers working in logic (the logical empiricists in
the 1920s and 1930s, especially Carnap and early Quine). Hilbert's
adoption of~$\lor$ in his \emph{Grundzüge der theoretischen Logic}
guaranteed its widespread use by mathematical logicians. The origins
of other logical symbols are also discussed.
\end{abstract}

\setlength{\epigraphwidth}{.7\textwidth}
\epigraph{%
No topic which we have discussed approaches closer to the
problem of a uniform and universal language in mathematics than does
the topic of symbolic logic. The problem of efficient and uniform
notations is perhaps the most serious one facing the mathematical
public. No group of workers has been more active in the endeavor to
find a solution of that problem than those who have busied themselves
with symbolic logic \dots\@ Each proposed a list of symbols, with the
hope, no doubt, that mathematicians in general would adopt them. That
expectation has not been realized.}{\citet[p.~314]{Cajori1929}}

\section{The mystery of `$\boldsymbol\vee$'}

The symbol `$\lor$' for inclusive disjunction is almost universally
accepted.  This contrasts with every other logical operator: negation
is symbolized by `$\lnot$', `$\pmnot$', `$\boldsymbol{-}$', or an
overline, conjunction is symbolized by `$\land$', `$\&$', and
`$\pmand$', and (material) implication is symbolized by
`$\boldsymbol{\rightarrow}$' and~`$\pmimp$'. So where do we get
`$\lor$' from and why have logicians more or less settled on this one
symbol? We answer this question here and show that we get the ubiquitous 
`$\lor$' for disjunction from Whitehead and Russell's 
pre-\textit{Principia} work in formal logic.

Our answer to this question differs from that given by earlier studies
of formal logic and its notations. The earliest use of `$\lor$' for
disjunction recorded in the encyclopedic list given by
\citet[307]{Cajori1929} is in Whitehead and Russell's 1910
\textit{Principia}. The Kneales say that ``the system [which includes
`$\lor$'] is that introduced by Peano in his \textit{Notations de
logique mathématique} of 1894, developed in the successive editions of
his \textit{Formulaire de mathématiques}, and then perfected by
Whitehead and Russell in their \textit{Principia Mathematica} of
1910'' \citep[p.~520]{Kneale1962}. In fact, however, Peano always used
`$\pmccup$' for inclusive disjunction and did not use `$\lor$' at all.
So one might suspect that the symbol for disjunction is one of the
improvements introduced by Whitehead and Russell in \emph{Principia
Mathematica}. But is that true? Why was `$\lor$' chosen specifically?
And why was it then universally adopted even as logicians introduced
alternatives for all of \emph{Principia}'s other logical symbols?

The story has two parts: The first part concerns the introduction of
`$\lor$' by Whitehead and Russell, who were working within the
notational and formal logical tradition of Peano.  We suggest here
that the initial introduction of `$\lor$' was motivated by a shared
commitment to two design principles for logical symbolisms
\citep{Schlimm2018,Schlimm2021,Toader2021}.\footnote{For an extended
discussion of these two design principles, see \citet{Schlimm2018} on
Frege and \citet{Schlimm2021} on Peano. We take the phrase ``design
principles'' from Schlimm. \citet{Toader2021} highlighted the second
principle.} The first principle is that symbols should be
\emph{unambiguous}---one shouldn't use the same symbol in different
meanings. The second principle is that symbols for analogous
notions---better, for notions satisfying analogous logical laws (say,
associativity)---should get \emph{similar} symbols: the symbols for
similar notions should have a similar, though not identical,
notation.\footnote{Peano puts this in terms of a \emph{principle of
permanence}: ``...when establishing a new system of notations, or a
new calculus, it is convenient to do it such that the new calculus be
similar \emph{as much as possible} with old calculi, so that the
student does not have to learn a whole new calculus, but only the
differences from the theory known to him'' (quoted from
\cite[85]{Toader2021}).}

As we suggest below, this first design principle led Peano to adopt a
new symbol, namely `$\pmccup$', to indicate a kind of logical union or
addition, and to avoid the then-prevalent symbol `$+$' for this use.
The point of this was to distinguish arithmetic and logical addition.
Russell followed Peano's practice here, as we will see
below.\footnote{Although Whitehead and Russell's pre-\textit{Principia} 
work in formal logic is our source for `$\lor$', we will sometimes speak 
of Whitehead and of Russell alone because the manuscripts or publications 
that we discuss are single-authored. Often we speak of Russell alone since 
Russell's manuscripts comprise most of the pre-\emph{Principia} materials 
that survive; indeed, the surviving Whitehead manuscript that we discuss 
below was kept in Russell's papers. Many of the pre-\emph{Principia} draft 
materials no longer survive. Yet in Russell's case (but sadly not in 
Whitehead's) the surviving materials provide a detailed record of his 
evolution both innotations and in design principles for notations. We do 
not think it is wrong to take the more extensive surviving textual record 
left behind by Russell as a record of Whitehead's evolving views on notation 
and design principles for notation pre-\textit{Principia} because Whitehead 
and Russell of course had to reach agreement on such points before
\textit{Principia} was published. See also the discussion in
Section~\ref{sec:white-russell}. Note that, as \citet[137--138]{Russell1948} 
recalls, Whitehead invented a great deal of the symbolism in \textit{Principia}. 
See also the remarks by \citet[\S V]{Urquhart1994} and \citet[\S XI]{Moore2014}.}
Subsequently, in the early logical investigations which ultimately led
to \emph{Principia}, Russell introduced the symbol `$\lor$' motivated
by that same design principle that Peano held: Russell wanted to
distinguish notationally between two kinds of logical addition,
namely, propositional disjunction (a.k.a. ``propositional addition'')
and class union (a.k.a. ``class addition''). In this, Russell
out-Peanoed Peano, who used `$\pmccup$' for both.\footnote{It should
be noted that neither Peano nor Russell followed this principle quite
perfectly. Peano used `$\fmimp$' for both implication and deductive
consequence \citep[\S3.2.1]{Schlimm2021}. And \emph{Principia} used
square dots both for conjunction and for scope indications.} But
Russell followed the second design principle, that analogous notions
should have similar symbols, in adopting `$\lor$' as a sharpening
of~`$\pmccup$' so as to stress the analogy between propositional sum
and class union. Russell did not think that this analogy was perfect.
He held that the formal laws for propositional addition and class
addition differed; hence he felt a need for a different, yet similar,
symbol.\footnote{\citet[\S\S2-3]{Bernstein1932} criticized the view
that these two notions are disanalogous in the way that
\emph{Principia} suggests in comments on $\pmast4\pmcdot78\pmcdot79$.}
But even where the formal laws are identical, Russell wanted to mark
the distinction between operations on different terms, as in the case
of class union and relation union. This pattern is noticeable in, for
instance, Russell's development of many symbols of his logic of
relations, which are usually distinguished only by an accenting mark,
suggesting that Russell tacitly followed the second design principle,
albeit a weaker version of it than Peano apparently endorsed.

All of that is the first part of the story. The second part of the
story is that of the adoption of $\lor$ in the period after
\emph{Principia}, when a number of logical symbolisms (and indeed,
logical methodologies) were in use simultaneously.  The main rivals of
\emph{Principia} were the algebraic approach of the
tradition of Boole, Peirce, and Schr\"oder, and the axiomatic logic of
Hilbert, with their respective notations. Despite being heavily
influenced by both \emph{Principia} and the algebraic notation of
Schr\"oder (which goes back to Peirce and the Boolean school), Hilbert
introduced different symbols for all the logical connectives---except
disjunction.

One upshot of the above findings is to uproot a misconception 
about the origins of the modern symbol for disjunction.
A common story in textbooks is that $\lor$ really is a version of the
letter~`v', the initial letter of \emph{vel}, the latin word for~``or''.
We wanted to find out whether this was a typographical accident, and
later writers simply thought the similarity between `$\lor$' and the
initial of \emph{vel} should be suggested as a mnemonic device to
students, or whether it actually played a role in the choice of the
symbol. Here we argue that although there is a single use of `v' as a
logical symbol abbreviating \emph{vel} in Leibniz, this played no 
known role in Peano's choice or in Whitehead and Russell's choice of 
notation for disjunction. However, this historical connection 
may, at least in part, have lead Hilbert to adopt it.

As usual in the history of mathematics, there is no simple answer to
the question of who deserves credit for~`$\lor$'.  In one sense, it
can be answered simply by recording the earliest use of the symbol.
But the earliest occurrence is often not the use to which the current
practice of the use of the symbol is historically connected. For
instance, Peirce in fact used~`$\lor$' for disjunction before Whitehead or Russell
did---but this use was independent and played no role in the
historical development of logical notation systems. (In this case we
know this for a fact because Peirce's relevant writings were not
published until long after the use of~`$\lor$' was already entrenched,
in 1933.) So it is only by tracing the adoption of Russell's symbol
that we can establish Russell as the originator of ``our''~`$\lor$'.

There is a related question, namely that of what should count as a
symbol for propositional disjunction in the first place.  Is Peano's
`$\pmccup$' a symbol for disjunction, or is it a symbol for something
else, perhaps an ambiguous combination of disjunction and union?  Who
invented disjunction in the sense we now use it---as a connective
between truth-apt expressions?  This is certainly a deeper and more
interesting question, but one we cannot answer here.  Our contribution
is more modest: to trace the history of the symbol `$\lor$' only, and
not the concept it expresses.\footnote{Our use of `the concept' in the
singular is a stylistic convenience: we do not mean to imply that
there is a unique concept throughout the history of logic that
correlates naturally to propositional disjunction and that is now
widely symbolized using `$\lor$'. For an in-depth treatment of the
concepts of disjunction, see \citet{Jennings1994}.} We do this not by
merely recording its first use, but by placing it in the context of
the development of logical symbolisms. Of course, even this is a much
larger context which we cannot do justice here---we focus only on the
parts that are relevant to~`$\lor$'.

Both of these considerations lead us to reject the claim that
Leibniz's use of~`v' should be considered the origin of~`$\lor$'. As
we will show, (a)~it had no known influence on the eventual
introduction and development of `$\pmccup$' and `$\lor$' by Peano, 
Whitehead, or Russell, and (b)~Leibniz's `v' was not a symbol for 
propositional disjunction.

We settle the questions of who invented the use of `$\lor$' for
disjunction (Whitehead or Russell) and whether this use of  `$\lor$'
was in fact introduced as an abbreviation for \emph{vel} (it was not). 
We suggest, however, that even though these questions may seem trivial,
they are not uninteresting, and we will mention, along the way, other
interesting aspects of the development of logical symbolism. Further,
the story serves as a case study in how the choice and adoption of
symbols and the organization of a symbolic system (a kind of
conceptual framework) can be determined by methodological commitments
but also by extrinsic factors such as the availability or cost of
printer's type or the influence that textbooks have on the adoption of
symbolisms. In this instance, in addition to the commitment to
unambiguous notations in Peano and Russell mentioned above, another
factor concerns the reasons that led Russell to choose
disjunction as a primitive of the logical system of \emph{Principia}.

\section{Leibniz}

Leibniz deserves some credit for the first use of `v' as a symbol for
`or' in some sense.\footnote{We make no claim here about Leibniz's
chronological precursors. We did, however, consult three medievalists
who specialize in history of logic about antecedent uses of $\lor$ for
inclusive disjunction, who did not recall any earlier uses.} But
Leibniz's use of~`v' is not for propositional disjunction, and occurs
in only one manuscript, his 1679 ``Matheseos universalis pars prior:
De terminis incomplexis'':
\begin{quote}
Also, in the way that $+$ is a conjunctive mark or sign of aggregation
corresponding to the word \emph{and}---as in $a + b$, that is, $a$ and
$b$ together---so a disjunctive mark or sign of alternation is also
given, which corresponds to the word \emph{or} [\emph{vel}], so $a
\mathbin{\leibnizv} b$ signifies $a$ or~$b$ to me. It also
has a use in calculus, for if we have $xx + ab = \overline{a + b}x$,
then $x$ will be equal to $a \mathbin{\leibnizv} b$,
whether it signifies $a$ or whether it signifies~$b$.  In this respect
it will have an ambiguous value.  For example, if $xx + 6 = 5x$,
$x$~can be~$2$, but $x$~can also be~$3$.  For if $x$ is~$2$, then from
$xx + 6 = 5x$, we get $4 + 6 = 10$; and if $x$ is~$3$, then from $xx =
6 = 5x$, we get $9 + 6 = 15$.  However, many values of this unknown
[variable] or of the present equation cannot be given roots, as will
be clear in the appropriate place below.

Hence, ambiguous signs also have a use, and this will be clear in the
appropriate place below, ambiguity being the font of irrationality in
calculus; and so when I write $x = 3 + \sqrt 4$, this can be
explicated by $3 + \sqrt 4$ or $3 + 2$ or~$5$, just as when I
write $3 - \sqrt 4$, it is explicated either by $3 - 2$ or~$1$;
thus, we get $x = 5 \mathbin{\leibnizv} 1$.  For in order
for us to remove the irrationality, let $x - 3 = \sqrt 4$; therefore
$xx - 6x + 9 = 4$, whether $xx - 6x + 5 = 0$ or $xx + 5 = 6x$, where
it is clear that it would satisfy~$5$ as much as~$1$.  For if $x$ were
given the value~$5$, we would get $25 + 5 = 30$; if $x$ is given the
value~$1$, we get $1 + 5 =
6$.\footnote{``\foreignlanguage{latin}{Quemadmodum etiam $+$ est nota
conjunctiva seu cumulationis et respondet \textgreek{τῷ} \emph{et}, ut
$a + b$ id est $a$ et~$b$ simul, ita datur quoque nota disjunctiva seu
alternationis quae respondet \textgreek{τῷ} \emph{vel}, sic $a
\mathbin{\leibnizv} b$ mihi significat $a$ vel~$b$. Idque
et in calculo usum habet, nam si sit $xx + ab = \overline{a + b}x$,
erit $x =  a \mathbin{\leibnizv} b$ seu $x$ significabit
vel $a$ vel~$b$ habebitque adeo valorem ambiguum. Ex.~causa si sit
$xx+6=5x$, potest $x$ esse~$2$, sed tamen potest etiam $x$ esse~$3$.
Nam si $x$ sit~$2$, tunc ex $xx+6=5x$ fiet $4+6=10$; et si $x$
sit~$3$, tunc ex $xx+6=5x$ lit $9+6=15$. Plures autem incognitae hujus
valores seu praesentis aequationis radices dari non possunt, ut suo
loco patebit.
\par
Hinc usum quoque habent signa ambigua, et suo loco patebit,
ambiguitatem in calculo esse fontem irrationalitatis; itaque cum
scribo $x = 3 + \sqrt{4}$, tunc id potest explicari tam per $3 +
\sqrt{4}$ seu $3+2$ seu~$5$, quam per $3-\sqrt{4}$ seu $3-2$ seu~$1$,
adeoque erit $x = 5 \mathbin{\leibnizv} 1$. Nam ut
tollamus irrationalitatem, sit $x-3=\sqrt 4$; ergo $xx-6x+9=4$ seu
$xx-6x+5=0$ seu $xx+5=6x$, ubi patet satisfacere tam $5$ quam~$4$. Nam
si $x$ valeat~$5$, fiet $25+5=30$; sin [sic] $x$ valeat~$4$, fit
$1+5=6$}'' (\cite[9v \& 9r]{Leibniz1679} and
\cite[57--58]{Leibniz1863}, translation courtesy of Jack Zupko).
Leibniz seems to be using the Greek article \textgreek{τῷ} to
indicate mention rather than use; Leibniz underlined the corresponding
occurrences of `et' and `vel' in the first sentence.}
\end{quote}

\begin{figure}
    \includegraphics[width=\textwidth]{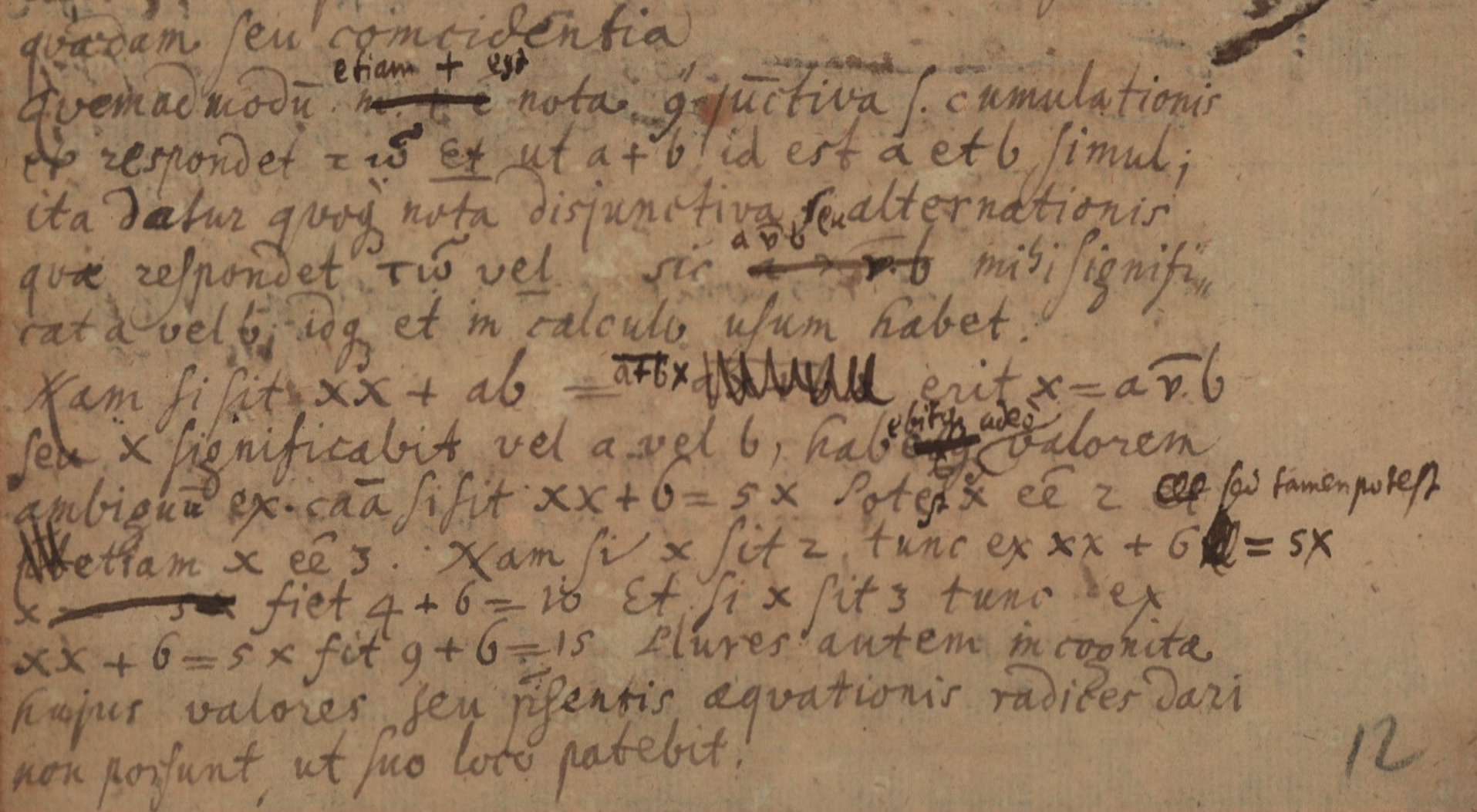}
    \caption{Use of $\mathbin{\leibnizv}$ in \citet{Leibniz1679}}
\end{figure}

This occurrence of `$\leibnizv$' is the first use of a `$\lor$'-like
symbol for `or' (in some sense) that we know of.\footnote{Why did he
put an inverse breve over~`v'? We conjecture that Leibniz wanted
to distinguish the constant near-nonce symbol `$\leibnizv$' from the
variable~`v' that he used frequently (e.g., \cite[57,
59]{Leibniz1679}). Another possibility is that it indicates an
abbreviation, like a scribal \emph{siglum}. However, the inverse breve
was not usually used for this purpose, and the standard abbreviation
for \emph{vel} was a lowercase `l' (ell), not~`v', combined with a vertical
stroke (\ipabar{l}{1ex}{1.1}{}{}) or macron (\= l).} However, Leibniz
does not seem to use the symbol here as a sign for disjunction of
(truth-apt) propositions as is usually the case today. Indeed, the
surrounding context does not include any uses of `$a$' or~`$b$' for
propositions. Rather, Leibniz seems to be using `$a$' and~`$b$' here
as standing for (truth-inapt) terms, and using `$a \mathbin{\leibnizv}
b$' to indicate \textit{ambiguity}, that is, one of two alternatives
without a specific one being indicated, as in phrases like ``either
chocolate or vanilla'' or ``7 or~8.''\footnote{We thank the reviewer for
helpfully suggesting that we put this Leibnizian use of~`$\leibnizv$
in terms of ambiguity.} That use of~`$\mathbin{\leibnizv}$' for
alternation between truth-inapt terms is consistent with the other two
uses `$a \mathbin{\leibnizv} b$' and `$5 \mathbin{\leibnizv} 1$',
where the symbol is flanked by (variables for) numbers and serves to
produce an ambiguously referring term.\footnote{See also
\citet[195]{Mugnai2018} on disjunction in Leibniz, and
\citet{Mugnai2010} more generally on logic and mathematics in the 17th
century.}

These three uses, in this one manuscript, seem to be the only uses of
the symbol~`$\mathbin{\leibnizv}$' in the entirety of his
\emph{\foreignlanguage{german}{Mathematische Schriften}}
\citep{Leibniz1863}. We
also found it nowhere else in Leibniz's writings. So this manuscript
may record the only uses of~`$\leibnizv$' for alternation
of (truth-inapt) terms in Leibniz's corpus. The manuscript was not
edited or published until~1863. It is perhaps unsurprising, then, that
Leibniz's use of~`$\leibnizv$' for disjunction (which, so
far as we can tell, was confined to this one manuscript) found no
imitators: even \citet{Arnauld1683}'s \emph{La logique ou l'art de
penser}, which was an especially influential book in the 1700s and in
the 1800s before Boole's work, omits `v' for disjunction. Logicians in
the 1700s and 1800s largely used symbols, when they did at
all, for the relations between terms in categorical logic and not
propositional disjunction; indeed, most logicians either did not
symbolize propositional disjunction and argued it was eliminable as a
form of reasoning, or else did not discuss it all; see
\citep[339]{Hailperin2004} for discussion.\footnote{\citet[\S
IV.4]{Whately1831}, for example, argues that all disjunctions of
categorical propositions are conditional in meaning due to the
equivalence of `$A$ is $B$ or $C$ is~$D$' with `if $A$ is not~$B$,
then $C$ is~$D$; and if `$C$ is not~$D$, then $A$ is~$B$'; see also
\citet[Chapter IV, p. 110]{Mill1843}.}

To qualify this slightly, a close textual reading of one or another
figure may extract, induce, or produce a propositional logic in some
sense.\footnote{To take one notable example, \citet{Lenzen2004} argues
that a system of propositional logic can be extracted from one of
Leibniz's categorical logics, and that many of its axioms are
so-formulated that they might reasonable be recognized by us today
as commonplace principles in logic textbooks treating the
propositional calculus.} Even so, our interest lies in the story of
how `$\lor$' became so widespread as a symbol for disjunction in
modern usage. It is our contention that Leibniz and other logicians
working primarily with categorical logics---between Leibniz and Boole,
we might say---did not inspire or cause the usage that we now have
today.

\section{The 19th Century}

In the course of mathematically analyzing logical principles using
algebraic equations, Boole used `$+$' as a symbol for something akin
to~`$\lor$'. In Boole's algebra, however, `$+$' is not a symbol for
disjunction in our sense, since the expressions that flank `$+$' are
not signs for propositions or statements, but for concept
expressions,\footnote{``Literal symbols, as $x$, $y$, \&c.,
representing things as subjects of our conceptions''
\citep[27]{Boole1854}. It should be noted that we do not mean to imply
that Boole's intended interpretation of his
symbols was fixed throughout his life. Indeed, \citet{Peckhaus1997}
and others have argued that there are significant differences between
\citet{Boole1847} and \citet{Boole1854} in their view of what terms
pick out (concepts versus extensions of concepts). These differences
are important, but incidental to our discussion since, either way, the
terms are not truth-apt (even though they can be used to induce or
extract a propositional logic, or can be applied in a manner that
mirrors propositional logics, much as one can do with Leibniz's
categorical logic.)} and so has rightly been thought of as an algebra
of classes (as we might put it), or as an algebra of concepts in
extension (as Boole seems to put it). 

Boole can appropriately use `$1$' for the universe, multiplication for
class intersection, and subtraction for complementation
\citep[20]{Boole1847}. Although Boole's operators `$+$' and~`$\cdot$'
applied primarily to classes, he also considered symbols~`$x$' for
propositions as expressing the ``cases'' where the proposition~$X$ in
question is true. The expression `$x+y$' then expresses the union of
the cases where $X$ and~$Y$ are true---assuming, however, that these
classes are disjoint. Claims about the propositions are made using
equations involving $1$ and~$0$ (the empty class). Multiplication is
re-interpreted to define (the algebraic analogue of) propositional
conjunction; similarly subtraction is re-interpreted as (the algebraic
analogue of) propositional negation. Then these operations
so-interpreted are similarly used to define (the algebraic analogues
of) propositional implication and propositional disjunction.  In
Boole's system, that either $X$ is true or $Y$ is true would be
expressed using the equation $x+y-xy=1$ \citep[51]{Boole1847}. 

Of course, Boole's choice of `$+$' for disjunction (in
\cite{Boole1847,Boole1854}) had a more immediate influence on notation
for disjunction in the 1800s.\footnote{In a 1848 manuscript ``The
Nature of Logic,'' Boole instead used `$+$' for (the algebraic
analogue of) propositional conjunction, adding parenthetically,
``which indeed was originally a contraction of the Latin \emph{et}''
\citep[6]{Boole2013}. In an 1856 manuscript ``On the Foundations of
the Mathematical Theory of Logic and on the Philosophical
Interpretation of Its Methods and Processes,'' he repeats this point
and suggests that the use of~`$+$' for \emph{et} and of~`$-$' for
\emph{minus} originates ``in medieval latin manuscripts''
\citep[87]{Boole2013}.} Subsequent choices of symbols for operations
on classes were influenced by Boole even when logicians departed from
Boole's philosophical program to varying degrees. For example,
\citet{Jevons1883} uses `$\cdot|\cdot$' for non-disjoint union, to
distinguish it from Boole's~`$+$'.\footnote{In his 1864 ``Pure
Logic,'' Jevons followed Boole in using `$+$' for disjoint union
\citep[80]{Jevons1864}. He changes the notation in \emph{The
Principles of Science} for a similar reason that leads Peano to adopt
$\pmccup$ (see below): ``[T]his sign [$+$] should not be employed
unless there exists exact analogy between mathematical addition and
logical alternation. We shall find that the analogy is imperfect, and
that there is such profound difference between logical and
mathematical terms as should prevent our uniting them by the same
symbol. Accordingly I have chosen a sign $\cdot|\cdot$, which seems
aptly to suggest whatever degree of analogy may exist without implying
more'' \citep[68]{Jevons1883}.} \citet[9]{Peirce1870} similarly uses
`$+$' for disjoint union and adopts `$\peirceplus$' for non-disjoint
union.\footnote{\citet[iii--iv]{Peirce1883} noted that the debate over
whether `$+$' should be reserved for disjoint union or not had been
decided: he notes that the majority of authors ``declared in favor of
using the sign of addition to unite different terms ino one aggregate,
whether they be mutually exclusive or not.'' Only Boole and
\citet{Venn1881} insisted on disjoint union.} In ``The Simplest
Mathematics'' \citeyearpar{Peirce1902}, which was not published until
1933, Peirce actually used `$\lor$' interchangeably with his
`$\peirceor$' symbol for inclusive or. By 1933, of course, the use of
`$\lor$' for disjunction was already entrenched. Other important
writers on the algebra of logic, such as
\citet{Schroeder1877,Schroeder1890}, MacColl \citep{McColl1877}, and
Ladd-Franklin \citep{Ladd1883}, all use `$+$' for union, and so does
\citet{Whitehead1898} originally. Even \citet{Nagy1890}, writing in
Italian and familiar with Peano's work, sticks to Boole's notation. 

Common to the followers of Boole is not just the choice of symbols,
but also the exclusive focus on equations. Using equations, they
express relations between class terms, and provide translations of
syllogistic reasoning into such equations.  However, MacColl
\citep{McColl1877}, \citet{Peirce1885}, Ladd-Franklin
\citep{Ladd1883}, and \citet{Schroder1891} took the further step of
showing that terms could translate propositionally valid inferential
relationships. They allowed operators like~`$+$' to connect not just
terms for classes or concepts but also terms for propositions or
statements, and developed systems of formulas, not just of equations.
In their work, the symbol~`$+$' serves as a genuine propositional
disjunction (in the inclusive sense).\footnote{On the development of
algebraic approaches to logic in the 19th century see
\citet{Peckhaus1997,Peckhaus2009}; on Schröder's logic, see
\citet{Dipert1991,Peckhaus2004}; on Peirce, Jevons, and Schröder, see
\citet{Dipert1978}; on Peirce and Peano's influence on Schröder, see
\citet{Peckhaus2014}; on MacColl, see
\citet{Peckhaus1999,Anellis2011}.}

Working alongside the algebraic logicians were Frege and Peano, who
both were important sources for Russell. Frege had his own
two-dimensional notation (first introduced in \textit{Begriffsschrift}
and later developed in \textit{Grundgesetze}). He was able to define
disjunction using the content strokes for implication and negation,
but did not introduce a separate symbol for it \citep[11]{Frege1879}.
Regarding Frege's notational influence on the story of~`$\lor$', the
record shows that Frege's actual notations and two-dimensional style
were rarely adopted (and sometimes were ignored or even maligned) by
others (with the exception of his popular turnstile `$\pmthm$'). But
Frege's notational influence on Russell was nonetheless felt: indeed,
Frege's views on \emph{design principles} for notation seems to have
been quite influential on Russell, even though his choice of actual
symbols was not. Russell's engagement with Frege's writings seems to  
have led Russell to later distinguish propositional sum and class sum,
for which purpose he introduced our~`$\lor$'---as we will see. First,
however, we consider the influence of Peano's notations and logical
framework on the development of \emph{Principia}.

\section{Peano}

Peano is responsible for much of the notations in \emph{Principia},
and Russell's notational imitations of Peano likely led him to
introduce `$\lor$' for disjunction. Indeed, Peano seems to be
responsible for suggesting to later readers that `$\lor$' was a
sensible (Latin-inspired) choice of symbol for disjunction. 

Peano, however, did not use~`$\lor$' (or the letter `v') for
disjunction or class union.  From his \emph{Calcolo geometrico}
onward, he instead used `$\pmccup$' for both (and `$\pmccap$' for
conjunction and class intersection).
\begin{quote}
    2. By the expression $A \pmccap B \pmccap C\pmccap \dots$, or
    $ABC\dots$, we mean the largest class contained in the classes
    $A$, $B$, $C$,~\dots{} or the class formed by all the entities
    which are at the same time in $A$ and $B$ and $C$, etc. The sign
    $\pmccap$ is read \emph{and}; the operation is indicated by the
    sign $\pmccap$ is logical \emph{conjunction}. We shall also call
    it \emph{logical multiplication}, and say that the classes $A, B,
    \dots$ are \emph{factors} of the \emph{product} $AB\dots$ 
    
    3. By the expression $A \pmccup B \pmccup C \dots$, we mean the
    smallest class which contains the classes $A$, $B$, $C$,~\dots{} or
    the class formed by all the entities which are in $A$ or $B$ or
    $C$, etc. The sign $\pmccup$ is read \emph{or}; the operation is
    indicated by the sign $\pmccup$ is logical \emph{disjunction}. We
    shall also call it \emph{logical addition}, and say that the
    classes $A$, $B$,~\dots{} are \emph{terms} of the \emph{sum} $A
    \pmccup B \pmccup \dots$ (\cite[1--2]{Peano1888};
    \cite[76]{Kennedy1973})
\end{quote}
Similarly, Peano  uses `$-A$' and `$\overline{A}$' for set
complementation and interprets it also as negation; he uses `$\bigcirc$'
for the empty class and interprets $A=\bigcirc$ as ``there is
no~$A$;'' and he uses `\rotatebox[origin=c]{30}{$\circlevertfill$}'
for the universal class and interprets $A \pmccup B =
\rotatebox[origin=c]{30}{$\circlevertfill$}$ \label{PeanoParagraph}
as ``everything is an $A$ or a $B$'' \citep[76--77]{Peano1888}. 

Peano was of course aware of the then-entrenched use of `$+$' and
`$\cdot$' or `$\times$' for both disjunction and conjunction and for
union
and intersection in the algebra of logic tradition. He deliberately
avoided using the same symbols in his own work, however:
\begin{quote}
    It seemed useful to substitute the symbols $\pmccap$, $\pmccup$, $- A$,
    $\bigcirc$, \rotatebox[origin=c]{30}{$\circlevertfill$} for the
    logical symbols $\times$, $+$, $A_\mathrm{i}$, $0$, $1$ used by
    Schr\"oder, in order to forestall any possible confusion between
    the symbols of logic and those of mathematics (a thing otherwise
    advised by Schr\"oder
    himself).\footnote{``\foreignlanguage{italian}Credetti utile di
    sostiture i segni $\pmccap$, $\pmccup$, $- A$, $\bigcirc$,
    \rotatebox[origin=c]{30}{$\circlevertfill$}, ai segni di logica
    $\times$, $+$, $A_\mathrm{i}$, $0$, $1$, usati dallo Schröder,
    affine d'imedire una possibile confusione fra i segni della logica
    e quelli della matematica (cosa del resto avvertita dallo
    Schr\"oder stesso)'' (\cite[x]{Peano1888}; translation from
    \cite{Kennedy1973}).}
\end{quote}

The symbols `$\pmccup$' and `$\pmccap$' are most likely adopted as
typographical variants of `$\grasssmile$' and~`$\grassfrown$'
introduced by \citet{Grassmann1844}.  Grassmann had used
`$\grassfrown$' as both a sign for an arbitrary operation
(\emph{Verkn\"upfung}) and as the sign for a kind of geometrical
product.  Leibniz had also sometimes used~`$\grassfrown$' as a symbol
for multiplication (e.g., \cite[5]{Leibniz1690}, reprinted in
\cite[15]{Leibniz1858b}, which Peano cites in \emph{Calcolo}).
Although Peano does not say, at least not in \emph{Calcolo}, these
uses may have suggested to him that `$\pmccap$' (or `$\grassfrown$')
would be suitable to express logical multiplication (conjunction). He
mentions Leibniz's use of `$\grassfrown$' for multiplication a decade
later in the \emph{Formulaire}:
\begin{quote}
    The sign $\grassfrown$ was adopted by Leibniz to indicate the
    arithmetical product. We adopt it for the logical product, since
    the arithmetical product is today indicated
    by~$\times$.\footnote{``\foreignlanguage{french}{Le signe $\grassfrown$
    a été adopté par Leibniz pour indiquer le produit arithmétique.
    Nous l'adoptons pour le produit logique; car le produit
    arithmétique est aujourd'hui indiqué par~$\times$}''
    \citep[32, our translation]{Peano1897}.}
\end{quote}

The propositional (re)interpretation of set theoretic operators is
purposeful: he extends (really, reaffirms) his propositional
interpretation of the set theoretic symbols in his treatment of
propositions in subsections 4 through 11 of \S1 on deductive logic, as
was common practice among (other than Frege) logicians in the 1800s.
Still, Peano explictly uses `$\pmccup$' for propositional disjunction,
e.g.:
\begin{quote}
    $\alpha \pmccup \beta$ expresses the condition that $\alpha$ is true or $\beta$ is true. \citep[82]{Peano1888}
\end{quote}
Peano's lower-case Greek letters may be (re)interpreted such that 
$\alpha \pmccup \beta$ is a union of classes wherever convenient. 
This sort of dual usage of symbols for classes or propositions, 
and for class operations or propositional connectives, allowed 
him to double-count proofs for analogous theorems.

Peano persists in this dual usage in \emph{The Principles of
Arithmetic} \citep{Peano1889}, where he used `$\pmccup$' as disjunction
and class union, and `$\pmimp$' for both subset and implication
\citep[105, 108]{Peano1889}. Such dual usage drew a lengthy
impeachment from \citet[242--247]{Frege1896}. But in fairness to
Peano, it leveraged Boole-style algebraic treatments in multiple
domains at once, and Peano had philosophical reasons to so-use
symbols.\footnote{Peano thought that it had been established by
``Boole, Schr\"oder, and others'' that deductive logic, like
arithmetic algebra, was part of algebra, that is, the ``calculus of
operations'' \citep[153]{Peano1891}. He even wanted to argue that
mathematical analysis produces results that mere logic would not
\citep[154, n.~3]{Peano1891}.} By his principle of permanence, since
the formal laws for propositional and class operators were identical
(in his view), Peano rightly used the same notation for these
different operators. This does not violate the first design principle
of symbolism prohibiting ambiguous symbols: any symbol in a theorem or
its proof must be unambiguous and the syntax is unimpeachable on this
score. The intended interpretation, on the other hand, may vary in
those limited cases such that the formal laws---even the theorems and
their proofs---do not differ. Indeed, Peano was generally very
thoughtful in his selection and use of symbols; see
\citet{Schlimm2021} for an insightful discussion. The ``design
principles'' Peano employed explain, for instance, why he chose pairs
of inverted symbols for dual operations like `$\pmccup$'
and~`$\pmccap$'.\footnote{Compare also the discussion of Frege's
symbolism in \citet{Schlimm2018}.}

Peano's habit of using one symbol for different notions---that is, for
allowing the intended interpretation to vary in order to economize the
amount of thought required to master the new symbols---is important to
note because it helps explain the genealogy of~`$\lor$'. Russell will
at first copy Peano's practice, although Russell was more
uncomfortable than Peano with using the same symbol under two distinct
intended interpretations. Even where the formal laws were identical,
as in the case of class notions and their relational analogues,
\emph{Principia} avoids using the exact same symbol for two different
notions. Rather, the distinction is marked in the syntax of the system
by using a different symbol entirely for different kinds of terms, as
we will see in the next section. Still, Russell did use similar
symbols for notions whose formal laws were close analogues,
distinguishing these different symbols often by a small change, like a
dot or an accent, in a Peano-like fashion. Thus Russell first used
`$\lor$' for a set-theoretic notion (relational product) and for
propositional disjunction, and then later uses it just for
propositional disjunction.\footnote{Peano would also write `$ab$' for
`$a \pmccap b$'. This notation is combined with the dot notations for
scope. Thus, \citet[3]{Peano1901} sometimes uses `$a \pmdot bc$' for
`$a(bc)$,'' which itself expands into ``$a \pmccap (b \pmccap c)$.''} What
Russell never does is to use the exact same symbol for distinct
intended interpretations: these are always marked in the formal
grammar by a different, even if highly similar, symbol. 

This is not to impeach Peano's practice. Peano may have felt justified
in his ambiguous use of `$\pmccup$' and~`$\pmccap$' by
the principle of permanence, according to which operations obeying the
same laws in different domains \emph{should} be symbolized by the same
notation (e.g., `$+$' and~`$\cdot$' applied to different sets of
numbers, say, integers and reals), and if not, they should be
symbolized differently. Since the laws governing class union and
intersection, and the laws for propositional disjunction and
conjunction are the same, this justifies the use of the same symbols
for both. By contrast, the arithmetical operations $+$ and~$\cdot$
obey different laws, e.g., idempotence fails (see \cite{Toader2021}
for discussion). Moreover, context disambiguates between the different
meanings of `$\pmccup$' and~`$\pmccap$'. By contrast, Peano's
class-theoretic construction of numbers in \emph{Arithmetices
principia} and later the \emph{Formulaire} has him using both `$+$'
and~`$\pmccup$' between class terms, and so using the same symbol
here would lead to ineliminable
ambiguity.\footnote{\citet[460]{Schroder1905} considers the
typographic contrast between Peano's logical and arithmetical
operators a benefit, as well as the horizontal symmetry between the
symbols for dual notions. He also suggests that one might best
remember that `$\grasssmile$' stands for disjunction by imagining the
symbol as a see-saw between two alternatives.}

We have already mentioned that Peano adapted his symbol~`$\pmccup$'
from Grassmann's~`$\grasssmile$'.  Is there any evidence that Peano was influenced in his
choice also by the similarity to the initial \emph{v} of \emph{vel},
or specifically by Leibniz's use of~`\leibnizv{}' for ``or''? Peano
was already aware of Leibniz's use of~`$\leibnizv$' as a symbol for
disjunction when he introduced his own: he cites the very page of
\emph{\foreignlanguage{german}{Mathematische Schriften}} where it
occurs in \emph{Calcolo} \citep[x]{Peano1888}. In his subsequent
writings on logic, he often mentions the fact
that Leibniz had used `v' this way (though Peano omits the accent),
and specifically that he did it because `v' is the initial letter of
\emph{vel}:\footnote{Peano deserves great credit for his keen eye: we
found only three occurrences of `$\leibnizv$' for disjunction in the
entire 393-page Volume 7 of Leibniz's mathematical writings (and we
found none elsewhere in the Gerhardt edition of Leibniz's mathematical
and philosophical writings). Leibniz is not explicit about why he
chose `$\leibnizv$', but Peano's guess is plausible and supported by
Leibniz's writing in Latin.}
\begin{quote}
Instead of $a\pmccup b$ Leibniz has $a \mathbin{\mathrm{u}} b$ (where
\emph{u} is the initial of
\emph{uel})\dots\footnote{``\foreignlanguage{italian}{Invece di $a
\pmccup b$ si scrisse da Leibniz $a \mathbin{\mathrm{u}} b$ (ove
\emph{u} `e l'iniziale di \emph{uel}), da Jevons $a \mathbin{.|.} b$,
e dal maggior numero di Autori $a + b$}'' (\cite[9]{Peano1891};
translation from \cite[158]{Kennedy1973}). Note that Peano uses the
variant spelling of \emph{vel} with an initial `u'---classical Latin
did not have separate letters for `U' and `V'. Presumably he did so
because `u' is more suggestive of Peano's chosen symbol~$\pmccup$. When
writing in Latin \citep{Peano1889} or later in his simplified
\emph{Latino sine flexione}---see \citep[Chapter 15]{Kennedy1980} for
discussion---Peano always spelled \emph{vel} with a~`v' whenever he
used the word.}

The sign $\circ$ corresponds to Latin \emph{aut}; the sign $\pmccup$ to
\emph{vel}.\footnote{``\foreignlanguage{french}{Le signe $\circ$
correspond au latin \emph{aut}; le signe $\pmccup$ à \emph{vel}}''
\citep[10]{Peano1894}.}

Leibniz, in his \emph{\foreignlanguage{german}{Schriften}}, vol.~VII
p.~57, indicates logical addition by `$a \mathbin{\mathrm{u}} b$'; the
letter u is the initial of the word
\emph{uel}.\footnote{``\foreignlanguage{french} Leibniz, dans ses
Œuvres, t. VII p. 57 indique l'addition logique par <<$a
\mathbin{\mathrm{u}} b$>>; la lettre u est l'initiale du mot
<<uel>>.'' (\cite[42]{Peano1897}, our translation).}
\end{quote}
However, Peano avoided a symbol reminiscent of `v for disjunction to
avoid confusion with other symbols: indeed, Peano already used
`\textsc{v}' for the universal class (\cite[139]{Peano1895};
\cite[13]{Peano1897}). Peano merely suggests that `$\pmccup$' for
disjunction can be thought of as a deformation of the `v' in
\emph{vel}:\footnote{A `v' symbol was also liable to be confused with
`\texttt{v}' for \emph{verum} used by \citet{Peirce1885} and Boole et
al.'s use of \emph{v} as a kind of existential quantifier.}
\begin{quote}
    $\pmast1\pmcdot2$ Let $p$ and $q$ be P [propositions]; $p \pmccup q$
    means ``at least one of these P is true''

    One may read the sign ``$\pmccup$'' as ``or''; this operation is called
    logical addition. [\dots]

    \emph{Leibniz} has indicated the operation~$\pmccup$ by the sign~$+$
    \citep[229]{Leibniz1890}; later (p.~237) by the same sign inside a
    circle. We cannot use a common sign for arithmetical and logical
    additions without producing ambiguities \dots 

    One may consider the sign $\pmccup$ as a deformation of \emph{v},
    the initial letter of ``vel,'' used by Leibniz as well for the
    same purpose. However, in the ``Arithmetices principia''
    \citep{Peano1889}, which contains the first theory rendered in
    symbols, I have chosen the shape of logical signs in such a way as
    to avoid any confusion.\footnote{``\foreignlanguage{french}{Soient
    $p$ et $q$ des P; $p \pmccup q$ signifie <<l'une, au moins, des
    ces P est vraie>>.\par On peut lire <<ou>> le signe $\pmccup$;
    cette operation s'appelle addition logique.\protect\par Ex. de la
    somme de Cls: \S Np P$2\pmdot1$ \dots; de P: \S$\times$ P$1\pmdot
    5$ \S Np P$1\pmdot2$ \dots\protect\par \emph{Leibniz} a indiqué
    l'opération~$\pmccup$ d'abord par le signe~$+$ (p.~229); ensuite
    (p.~237) par la même signe dans un cercle. Nous ne pouvons pas
    représenter par un même signe les additions logique et
    arithmétique, sans produire des ambiguïtés. Voir
    \S$+$P$5$.\protect\par On peut considérer le signe $\pmccup$ comme
    une déformation de \emph{v}, lettre initiale de <<vel>>, aussi
    employée par Leibniz pour le même but. Mais dans les
    <<Arithmetices principia, a.1889>>, qui contient la première
    théorie réduite en symboles, j'ai fixé la forme des signes de
    Logique de façon à éviter toute confusion.}''
    (\cite[16]{Peano1899}, our translation).\protect\par Note that
    here Peano does use `v' and the spelling \emph{vel}. Leibniz's use
    of `$+$' referred to here by Peano is for class union not
    propositional disjunction.}
\end{quote}

So even though Peano often suggested a connection between the symbol
for disjunction and the Latin \emph{vel}, it is unlikely that his
choice of~`$\pmccup$' for ``or'' was originally motivated by this
connection. Russell likely got the idea of using `$\lor$' for
disjunction from Peano, whom he read closely.\footnote{``In July 1900,
there was an International Congress of Philosophy in Paris\dots\@ The
Congress was a turning point in my intellectual life, because I there
met Peano. I already knew him by name and had seen some of his work,
but had not taken the trouble to master his notation. In discussions
at the Congress I observed that he was always more precise than anyone
else, and that he invariably got the better of any argument upon which
he embarked. As the days went by, I decided that this must be owing to
his mathematical logic. I therefore got him to give me all his works,
and as soon as the Congress was over I retired to Fernhurst to study
quietly every word written by him and his disciples. It became clear
to me that his notation afforded an instrument of logical analysis
such as I had been seeking for years, and that by studying him I was
acquiring a new and powerful technique for the work that I had long
wanted to do'' \citep[217--218]{Russell1967}.} This is not terribly
suprising: indeed, it is hard to overstate Peano's influence on
\emph{Principia}'s notations. Even \citet{Whitehead1898} followed
Boole and symbolized disjunction and conjunction as addition and
multiplication before encountering Peano's work. After becoming
thoroughly familiar with Peano's work through Russell,
\citet{Whitehead1902} also adopted Peano's `$\pmccap$' and~`$\pmccup$'.
\citet[vii]{WhiteheadRussell1910} fully acknowledge Peano's notational
impact: ``In the matter of notation, we have as far as possible
followed Peano\dots'' 

Yet in using `$\lor$' for disjunction, \emph{Principia} departs from
Peano's notations. As we will see in the next section, Whitehead and
Russell are jointly responsible for this. The earliest surviving use
of~`$\lor$' for disjunction we could find is in a Whitehead
manuscript, but Russell may well be the one who introduced the symbol.
He first separated the notations for class union and propositional disjunction.
Russell also used `$\lor$' for disjunction in solo published work
some years before the co-authored \emph{Principia} was published. The
next section supplies the evidence for this summary.

\section{Whitehead and Russell}\label{sec:white-russell}

We give Whitehead and Russell joint 
credit for first using `$\lor$' for disjunction. The relevant surviving 
manuscripts from Whitehead and Russell were written in close proximity, 
but the relative dating is uncertain. Moreover, many manuscripts from 
the time of the writing of \textit{Principia} are lost. In any case,
one of Whitehead and Russell, or both jointly, first used `$\lor$' 
for disjunction. Their use of this notation, through the influence of
\textit{Principia}, is the reason for today's wide acceptance of it
(as we argue in the next section). 

After framing the timeline for our discussion, we trace the genealogy
of~`$\lor$' in the surviving manuscripts. We will start from what is
definitely known and conclude this section with what is merely
probable. In our genealogy of Whitehead and Russell's design principles 
for notation we also focus on Russell's publications and manuscripts 
because few of Whitehead's pre-\emph{Principia} manuscripts from this 
period survive.\footnote{Of course, Whitehead no less than Russell was 
responsible for much of the notation in \emph{Principia} and both authors 
gave notation a great deal of thought. We also do not mean to suggest that 
their thoughts on notations and design principles for notations were
always identical 
in the pre-\emph{Principia} years; see \citep[137--138]{Russell1948} for 
some indication of notational divergences between the two. See also the 
letters quoted by \citet[\S V]{Urquhart1994} and \citet[\S XI]{Moore2014}.} 
Our aim is to correctly trace the evolution of \emph{Principia}'s notations, and the design principles for
choosing them. Crediting either Russell or Whitehead with any specific
decision is difficult for the above reasons. Thus, even where we
consider a single-authored manuscript or publication and write of,
e.g., as Russell introducing a symbol, we do not mean to imply that
Russell alone deserves the credit. We have no doubt that Whitehead and
Russell agreed on the design principles for notations as they are
manifested in \textit{Principia}. 

Here is the time frame for the remainder of this section: in his 1901
``The Logic of Relations'' \citeyearpar{Russell1901a}, Russell took implication as a primitive
symbol. In Whitehead and Russell's 1910 \emph{Principia},
`$\lor$'~for disjunction is now primitive
\citep[93]{WhiteheadRussell1910}. This raises two questions about
Russell's development: (1)~When and why was disjunction taken as
primitive? (2)~When and why was `$\lor$' introduced for disjunction?

We first take up~(1). To answer the ``why'' part: disjunction is
taken as primitive because the resulting system is simpler in that it
has fewer primitive propositions.\footnote{As to what ``simpler''
means here: \citet[82--83]{BellucciMoktefiPietarinen2018} have
discussed the various notions of simplicity at play in discussions of
Peano, Frege, and Peirce in their work on logic: in one sense of
``simpler'' (as applied to logical systems), a simpler system has
fewer total primitive symbols. Another sense of ``simpler''  is
that there are fewer primitive constants involved in the expression of
some formulas. Whitehead and Russell in \emph{Principia} are pursuing
simplicity in the sense of reducing the total number of primitive
symbols and primitive principles, and so are prioritizing simplicity
in the first sense. Note that unanalyzability is not required for
simplicity in Whitehead and Russell's sense, for
\citet[vi]{WhiteheadRussell1910} openly demur on whether ``the
analysis could not have been carried farther: we have no reason to
suppose that it is impossible to find simpler ideas and axioms by
means of which those with which we start could be defined and
demonstrated.'' So it cannot be that their sense of simplicity here is
tied to being unanalyzable. See also the quotes discussed in this
section, and especially the joy with which
\citet[xiii--xiv]{WhiteheadRussell1925} welcomed the analysis of the
1910 set of primitives by \citet{Sheffer1913} and \citet{Nicod1917}.}
In \emph{Principia}, \citet[130]{WhiteheadRussell1910}
comment on $\pmast5\pmcdot25$, which is $\pmthm \pmdottt p \pmor q
\pmdot \pmiff \pmdott p \pmimp q \pmdot \pmimp \pmdot q$:
\begin{quote}
    From $\pmast5\pmcdot25$ it appears that we might have taken
    implication, instead of disjunction, as a primitive idea, and have
    defined ``$p \pmor q$'' as meaning ``$p \pmimp q \pmdot \pmimp
    \pmdot q$.'' This course, however, requires more primitive
    propositions than are required by the method we have adopted.
\end{quote}
If we set aside inference rules and consider just propositional axioms
or axiom schemata, the texts bear out this claim. ``The Theory of
Implication'' \citep{Russell1906} gives seven axioms for propositional
logic; \emph{Principia} instead has five.\footnote{However, a judgment
as to overall simplicity somewhat depends on how we count. If we are
counting by \emph{primitive propositions}---including rules of
inference like \emph{modus ponens}, well-formedness rules, and uniform
substitution rules---then both systems have ten starred primitive
propositions. Hence, \citet[224]{Linsky2016} has a list of nine
primitive propositions in the 1906 system (by counting \emph{modus
ponens} and merging the two rules for uniform substitution);
\citet[226--227]{Linsky2016} also has a list of seven total primitive
propositions in \emph{Principia}, which omits the rules of
well-formedness $\pmast1\pmcdot7\pmcdot71\pmcdot72$.}

The ``when'' part of question (1) can be answered using the
Couturat-Russell letters translated in the headnote to
\citep{Russell1906}. The first known evidence of the switch from
`$\pmimp$' being primitive to `$\lor$' being primitive is from 21 August
1906.

Russell wrote ``The Theory of Implication'' in the summer of 1905, and
it was published in the \emph{American Journal of Mathematics} in
April~1906. Again, Russell adopts (material) implication as a
primitive therein. Russell justifies `$\pmimp$' being primitive
explicitly on much the same grounds as he later justifies replacing
`$\pmimp$' with~`$\lor$' as primitive: ``it avoids hypotheses which
are otherwise necessary'' \citep[162]{Russell1906}. What hypotheses
does Russell avoid in the 1906 paper? As a footnote makes clear
\citep[162]{Russell1906}, these hypotheses are to the effect
that for any $p$ and $q$, if $p \pmimp q$, then $p$ and~$q$ are
propositions, and so are truth-apt. These hypotheses are adopted in
\emph{The Principles of Mathematics} \citep{Russell1903} to prevent
the relation of material implication from obtaining between any
non-propositional (truth-inapt) entities.

Recall that Russell's ``calculus of propositions,'' also known as the
\emph{substitional theory}, allowed for variables to range over
everything without restriction \citep[224]{Linsky2016}.
\citet[558]{Landini1996}  calls this ``the doctrine of the
unrestricted variable'' because on this theory variables range over
all entities whatever (rather than, say, different sorts of variables
ranging over objects and predicates, as is now common in presentations
of second-order logic). This was only possible because the notions of
`$\pmimp$' in \citep{Russell1903,Russell1906} were quite different
from the modern one and that of the 1910 \emph{Principia}. On the
earlier substitutional theory, the sign~`$\pmimp$' is a relation sign:
it connects terms to make a term, so that ``Socrates $\pmimp$
Socrates'' was well-formed in Russell's early systems. The modern sign
for implication `$\rightarrow$,' like the symbol `$\pmimp$' in
\emph{Principia}, behaves very differently: it connects well-formed
formulas to make a well-formed formula, so that ``Socrates
$\rightarrow$ Socrates'' is syntactically ill-formed (much like 
``\emph{a} $\rightarrow$ \emph{a}'' is ungrammatical in most modern systems). 

Of course, the primitive postulates in \citet[\S18]{Russell1903} precluded 
any material implication relations obtaining between non-propositional
entities even though such strings could be syntactically well-formed:
one does not want to allow that ``Socrates implies Socrates'' is a 
truth of one's formal system (and this is not a theorem in either 
early system). How did Russell prevent, say, ``Socrates $\pmimp$ Socrates'' 
from being even vacuously true in the 1903 system while preserving the 
doctrine of the unrestricted variable under which any term is substitutable
for a variable? Russell built into some primitive propositions of the 1903 
system the antecedent condition that $p$ is a proposition. For example, 
in place of the axiom $p \pmimp p$, Russell had $p \pmimp q \pmdot
\pmimp \pmdot p \pmimp p$, which on Russell's intended interpretation is 
tantamount to adding the condition ``whatever implies anything is a
proposition'' \citep[\S18]{Russell1903}. 

In ``Theory of implication,'' \citep[$\pmast1\pmcdot2$]{Russell1906} 
instead interprets `$\pmimp$' in such a way that it may be just vacuously 
false whenever truth-inapt terms flank it, making the further hypothesis 
that whatever implies anything is a proposition redundant.\footnote{\citet[23]{Russell1906}
credits Frege's \emph{Grundgesetze}, Volume I with this shift in
interpretation of~`$\pmimp$.'} Russell's intended interpretation 
of the 1906 system effects some simplification of the primitive 
propositions in the 1906 system as compared with the 1903 one.\footnote{Note 
that the number of primitive propositions---again, depending on whether one 
counts axioms and rules of inference---is still ten in the 1903 system.}

At any rate, in 1906 Russell still takes `$\pmimp$' as primitive. 
Russell finally introduces `$\lor$' for disjunction in~$\pmast4\cdot11$
\citep[176]{Russell1906}, which defines $p \pmor q$ as $\pmnot p
\pmimp q$. The first passage of the section makes clear that `$\lor$'
means disjunction (``propositional sum''):
\begin{quote}
    In this section we shall be concerned with the fundamental
    properties of the \emph{propositional product} and the
    \emph{propositional sum} of two entities $p$ and $q$\dots\@ The
    \emph{propositional sum} of $p$ and $q$ is practically ``either
    $p$~is true or $q$~is true.'' We avoid a new primitive idea by
    taking as the propositional sum $\pmnot p \pmimp q$, i.e., ``if
    $p$~is not true, then $q$~is true.'' There are some advantages in
    taking the propositional sum as a primitive idea instead of `$p
    \pmimp q$'; we then define `$p \pmimp q$' as the propositional sum
    of $\pmnot p$ and~$q$. The choice is a matter of taste.
    \citep[175]{Russell1906}
\end{quote}
As this quotation shows, Russell was obviously not wedded to one or
another choice of primitive in his 1906 paper and allows that either
one can define the other. This is a significant shift from his view in
\emph{Principles}, where he explicitly argues that implication cannot
be defined, even in terms of equivalent formulas using other
connectives \citep[\S16]{Russell1903}.

Shortly after publication of ``Theory of implication,'' Russell found a reason to
prefer disjunction as a primitive notion over implication: it reduced
the number of axioms needed for propositional logic. Couturat and
Russell exchanged letters about Russell's 1906 paper. On 21 August
1906, Russell wrote to Couturat to announce his discovery that taking
`$\pmor$' as primitive simplifies the development of his logic:
\begin{quote}
    One can take as a primitive idea $p \pmor q$, while defining 
    \[ p \pmimp q \pmdot \pmiddf \pmdot \pmnot p \pmor q \] In this
    way two primitive propositions are made superfluous.
    \citep[\#216]{Schmid2001}
\end{quote}

This answers question (1): by 21 August 1906, Russell took disjunction
as primitive instead of material implication, and he did this because
it was logically simpler, that is, it reduced the number of axioms
needed in propositional logic.

Now for question (2): when and why did Whitehead and/or Russell first 
choose the symbol `$\lor$' for disjunction? Unlike Peano, Whitehead and 
Russell never explicitly explained this choice of notation. So we will 
have to reconstruct their reasons.

The two natural conjectures are that (A) Whitehead and Russell 
adopted~`$\lor$' because Leibniz used~`$\leibnizv$' and (B) Whitehead and 
Russell adopted~`$\lor$' as a modification of the already available 
notation~`$\pmccup$.' Given the available evidence, we endorse~(B). 
First let us consider why (A)~is improbable. 

If one supposes that Whitehead and Russell got `$\lor$' from Leibniz, 
the next question to answer is, `How did Whitehead and Russell come 
across it?' Here there are two natural possibilities: (A1) Whitehead 
and Russell found a use of~`$\leibnizv$' in Leibniz independently, or 
(A2) Whitehead and Russell found out about this possibility through 
Peano, who, as we saw, did notice Leibniz's use of~`$\leibnizv$'.

Against (A2), there is no textual evidence, for instance, in
\citet{Russell1900}, that Whitehead or Russell encountered, on their own, 
Leibniz's sole manuscript in which `$\leibnizv$' is used for disjunction.
Indeed, there is no evidence that Whitehead or Russell ever owned 
\emph{Leibnizens mathematische Schriften}. No notes on these volumes occur 
in Russell's Leibniz notebooks \citep{ArthurGriffin2017}. \emph{Leibnizens
mathematische Schriften} are also not in Russell's library, which is
held at the Bertrand Russell Research Centre. By contrast, Russell's
library has all seven volumes of the \emph{Philosophischen Schriften},
which he annotated heavily; these annotations were published by
\citet{ArthurGalaugherGriffin2017}. Furthermore, Russell kept a list
of books he was reading from 1891--1902 \citep{Russell1891}. None of
\emph{Leibnizens mathematische Schriften} appear on this list of books
\citep[165-166]{OBriant1979}. Finally, Russell's citations of
\emph{Leibnizens mathematische Schriften} are sparse enough to have
been pulled from other authors who quoted from them
directly.\footnote{Russell's appendix to the index of Leibniz quotes
includes just 18 citations to the entire \emph{Mathematische
Schriften}, and not even 10 to any single volume in that series. By
contrast, there are 492 citations to the \emph{Philosophischen
Schriften} series, and, to each of the seven volumes, 15, 143, 39, 63,
108, 62, and 62 citations, respectively.} So it is probable that
Russell never saw the manuscript in question, or even the volume in
which it was available.\footnote{Russell's Leibniz book has two
citations to \emph{Leibnizens mathematische Schriften}, Volume VII,
but none to the relevant manuscript \citep[109, n.~1,
247]{Russell1900}. And these two quotes were probably pulled from
\citet[64-65]{Cohen1883}, whom Russell cites in the footnote in
question. Similarly, in \emph{Principia}'s introduction and treatment
of propositional logic ($\pmast1$--$\pmast5$), Volume VII is cited
exactly once, in a footnote to $\pmast3\pmcdot47$, but again, the
relevant manuscript is not cited. This citation in \emph{Principia} is
also identical to an earlier citation; see \citep[182]{Russell1906}.}
Similar reasoning applies to Whitehead, whose 1898 \emph{A Treatise on 
Universal Algebra} cites Leibniz (as ``Leibnitz'') just once, discussing 
a theorem about resultant forces. And Whitehead reports, ``My knowledge 
of Leibniz's investigations [when I wrote the \emph{Treatise}] was 
entirely based on L. Couturat's book, \emph{La Logique de Leibniz}, 
published in 1901.'' \citep[10]{Whitehead1951} In summary, it is consistent 
with the available evidence that neither Whitehead nor Russell ever 
directly saw Leibniz's use of `$\leibnizv$'. 

It is of course true that Whitehead and Russell had some first-hand 
acquaintance with Latin and Greek thanks to their pre-college education. 
Russell was educated in Latin and Greek, in addition to being fluent in 
French and German. Might it be possible that Russell wanted to associate 
propositional disjunction with `\emph{vel}' even if he never saw the 
occurrence in Leibniz? We find this quite unlikely because Russell 
despised Latin (and Greek), writing in his \citep[36]{Russell1967}, 
``I hated Latin and Greek, and thought it merely foolish to learn a 
language that nobody speaks.'' (Happily, he did not feel similarly about 
logical languages.) Moreover, as \citet[105]{Lenz1987} notes, ``The learning 
of Latin and Greek played only a small role in [Russell's] tutelage.'' Even if
Peano's remarks about Leibniz may have caused Russell to think
of~`$\lor$' when deciding on a symbol to represent disjunction, there
is no reason to think that a connection to \textit{vel} motivated
Russell's choice of~`$\lor$', and some reason to think that Latin
associations would have been, if anything, wholly repugnant to him. 
Whitehead also studied Latin and Greek early on, but there is little 
evidence to suggest this connection with \emph{vel} would have occurred 
to Whitehead (or have felt apt to impress the new symbol's meaning upon 
the mathematical community). Indeed, Whitehead's education early on let 
Latin take a backseat to mathematics: ``I was excused in the composition 
of Latin Verse and the reading of some Latin poetry, in order to give 
more time for mathematics'' \citep[6]{Whitehead1951}.

On the other hand, in support of~(A2), Whitehead and Russell were close 
readers of Peano's works, including the \emph{Formulaire}, throughout the 
period from 1900 to~1903. In multiple works, Peano explicitly noted Leibniz's
use of~`$\leibnizv$'. So if one is to accept the view that `$\lor$'
comes to us from Leibniz, Peano is the most probable imtermediary
source of Whitehead and Russell's using~`$\lor$' for disjunction, even if Peano
never adopted Whitehead and Russell's notation for disjunction. But against this
proposal, there was a years-long delay between Whitehead and Russell's close study
of Peano's works and their first uses of~`$\lor$.' Russell studied
Peano's works extensively in September 1900, just after the July 1900
Congress \citep[363]{Russell1891}. Whitehead references the first three volumes 
of the \emph{Formulaire} and praises the ``admirable work'' of Peano's school 
in a footnote to his April 1901 ``Memoir on the Algebra of Symbolic Logic'' \citep[140]{Whitehead1901}. 
However, as we will see, Russell first used~`$\lor$' in his spring 1902 ``On Likeness'' manuscript, and 
Whitehead first uses it in his contemporaneous ``The Logic of Propositional Functions'' manuscript.
So neither Whitehead nor Russell used this notation until many months
after their deep study of Peano. 
This suggests it was not Leibniz's notational choices that led Whitehead and Russell
to adopt~`$\lor$' for disjunction, but some other development. We
describe just such a development, consistent with hypothesis~(B), that
explains why Whitehead and Russell modified `$\pmccup$' to use~`$\lor$' for
disjunction. Accordingly, we reject hypothesis~(A).

We saw above that Peano used `$\pmccup$' for both disjunction and class
union---not ambiguously in a problematic sense, but so that one could
follow a different intended interpretation of the symbol and retain
all the formal laws governing~`$\pmccup$' and the proofs involving
occurrences of~`$\pmccup$'. This accords with the ``design principle''
for notations, the principle of permanence, discussed above. Russell
seems to have accepted a similar principle, but not quite the same
one. Indeed, even where the formal laws are identical,
\emph{Principia} distinguishes between different kinds of terms. So in
\emph{Principia} we find `$\alpha \pmccup \beta$' for class union in
$\pmast 22$ and `$R \pmrcup S$' for relational union (i.e., a union of
ordered couples) in~$\pmast 23$. This is despite the fact that ``the
definitions and propositions of this number [$\pmast 23$] are to be
exact analogues of those of~$\pmast 22$''
\citep[226]{WhiteheadRussell1910}. 

So Russell cannot have held that different notions should get the same
notation even when the formal laws are identical. In rejecting Peano's
design principle here, Russell was likely further influenced by
Frege's sharp criticisms of Peano's conceptual notations.\footnote{We
say ``further influenced'' by Frege because Russell was already headed
in the direction of separating notations for different notions. By
spring 1902 Russell separated `$\lor$' for class union and `$\cup$'
for propositional disjunction \textit{before} engaging deeply with
Frege's writings in the summer of 1902.} Indeed, as we will see below,
the very manuscript wherein Russell first introduces `$\lor$' for
disjunction is written after Russell's engagement with Frege's
writings (although Russell seems to have been moving in this direction
already; see \citealt[xvii]{Urquhart1994}). Furthermore, Frege
explicitly criticizes Peano's alleged ``twofold use'' of symbols for
propositional and class notions.\footnote{``Not least disturbing, it
seems to me, is its falling apart into the calculus of classes and the
calculus of judgments, as it is customary to put it. And this
separation is already less marked in Mr.~Peano's work''
\citep[242]{Frege1984}.} So if Russell introduced `$\lor$' to
distinguish propositional sum and class sum, as we think he did, the
fact that Russell had just taken extensive notes on Frege's writings
(and Schr\"oder's writings---Schr\"oder likewise criticized Peano for
alleged dual uses of notation for propositional and class notions), it
is no coincidence \citep{Linsky2004}.\footnote{Nor is it a coincidence
that Russell, having just finished his detailed study of Frege (and of
Schr\"oder), writes in \textit{Principles} a Fregean (Schr\"oderian?)
broadside against Peano for not using separate symbols for
propositional and class notions: ``The subject of Symbolic Logic
consists of three parts, the calculus of propositions, the calculus of
classes, and the calculus of relations. Between the first two, there
is, within limits, a certain parallelism, which arises as follows: In
any symbolic expression, the letters may be interpreted as classes or
as propositions, and the relation of inclusion in the one case may be
replaced by that of formal implication in the other...A great deal has
been made of this duality, and in the later editions of the
\emph{Formulaire}, Peano appears to have sacrificed logical precision
to its preservation'' \citep[\S13]{Russell1903}.} All this supports
conjecture~(B): Russell introduced~`$\lor$' to separate notations for
two different notions---propositional sum and class sum---and was
further influenced by Frege's writings to do this.

Still, `$\lor$' did not appear \textit{ex nihilo}. Rather, this symbol 
was a modification of the symbol~`$\pmccup$'. Russell's apparent 
practice, even after engaging with Frege's work, was that \emph{similar} 
notions should have \emph{similar} notations. Admittedly, Russell never 
explicitly lays this down as a design principle for symbolism---he is 
far less explicit than Peano in this regard---but we know that Russell 
was, as he described himself, ``fussy'' about even minute aspects of how
his writing appeared, such as where punctuation and paragraph breaks
occurred, whether a colon or semicolon should be used, and formatting
issues like italics \citep{Blackwell1983}. Needless to say, Russell
was a consciously conscientious designer of notations, even if he
never made his design principles for notations explicit.

That Russell embraced this design principle can be inferred from his 
notational design principles of (1)~using similar symbols for similar 
notions (ones obeying analogous laws), and of (2)~using distinct (even 
if similar) notations for distinct notions. There are many examples in 
\emph{Principia} supporting that these two design principles were 
Russell's apparent (if implicit) practice. For example, we find `$\pmimp$' 
and~`$\pmcinc$' used in \emph{Principia} for implication and class
inclusion,
respectively. We find `2' and~`$\dot{2}$' for the cardinal number~2 and 
the ordinal number~2, respectively. The formal laws governing these pairs 
of notions differ; nonetheless, \emph{Principia} uses similar notations
for them because they are analogous. Crucially, \emph{Principia} never
uses the same symbol for distinct notions whose formal laws differ. 

This fact helps us answer the `why' part of question~(2): we maintain
that Russell chose `$\lor$' for disjunction because he was keen
to separate symbols for propositional and class notions because he
held that their formal laws differed. In this Russell was unlike Peano
who deliberately wed them together, since he held that the formal laws
that each obeyed were the same. Russell's choice of `$\lor$' thus
mimics Peano's practice of using some symbols dually for propositional
connectives and class operators: Russell wants to analogize the
notions expressed by `$\pmccup$' and `$\lor$'. But because the analogy
was imperfect, in Russell's view, Russell chose a new but similar
symbol for propositional disjunction, namely, a sharpened,
pointed~`$\pmccup$'.

Why did Russell think that `$\lor$' and `$\pmccup$' obeyed distinct
formal laws? In ``The Theory of Implication''
\citep[190]{Russell1906}, remarks:
\begin{quote}
    The analogues, for classes, of $\star 5\pmcdot76 \pmcdot 79$, are false.
    Take, \emph{e.g.}, $\star 5\pmcdot 78$, and put $p = $ English people, 
    $q = $ men, $r = $ women. Then $p$ is contained in $q$ or~$r$, but is 
    not contained in~$q$ and is not contained in~$r$.
\end{quote}
This exact remark is repeated in \emph{Principia} following 
$\pmast 4\pmcdot 78 \pmcdot 79$. And Russell's distinction of the 
propositional and class calculus goes back to \emph{Principles}. In 
\S25~of \emph{Principles} \citep{Russell1903}, titled ``Relation [of
the calculus of 
classes] to propositional calculus,'' Russell gives an example of a formula 
that ``can only be truly interpreted in the propositional calculus: in the 
class-calculus it is false.'' Russell then says this about disjunction:
\begin{quote}
    The disjunction is what I shall call a \emph{variable} disjunction, as 
    opposed to a constant one: that is, in some cases one alternative is 
    true, in others the other, whereas in a constant disjunction there is 
    one of the alternatives (though it is not stated which) that is always 
    true. Wherever disjunctions occur in regard to propositional functions, 
    they will only be transformable into statements in the class-calculus in 
    cases where the disjunction is constant.
\end{quote}
Though the vocabulary is unfamiliar, the underlying point is familiar
to modern readers, who might put it as follows: a disjunction of $p
\lor q$ might be such that $p$ is always true, as when $p$ is replaced
by `$0=0$,' or might be such that one of $p$ or $q$ is true, though
neither is always true, as when we have $p$ replaced by `$a$ is human'
or `$a$ is non-human.' The relevance of this passage is that it shows
Russell holds that the formal laws in the calculus of classes and
those in the calculus of propositions differ. Accordingly, as early as
his 1902 manuscripts where `$\lor$' is introduced, Russell wants a new
symbol to distinguish disjunction and class union.

The hypothesis that Russell chose `$\lor$' to logically distinguish 
disjunction as obeying different formal laws, while still preserving its 
analogy with~`$\pmccup$', is further supported by his letter to Frank Morley, 
the editor of the \emph{American Journal of Mathematics}, concerning the 
publication of his 1906 ``The Theory of Implication.'' The symbol~`$\lor$' 
and others were so unfamiliar to Morley that he wrote to Russell
about it and noting, perhaps ominously, that authors occasionally had to 
bear ``any \emph{heavy} expense for new type.'' Russell wrote back:
\begin{quote}
    Thank you for your letter. I do not think that any new type ought
    to be needed, seeing that you have already printed Whitehead's
    memoirs. I suppose the sign~`$\pmor$' is doubtful. It is very
    desirable that it should be pointed, not round; but otherwise it
    doesn't matter much whether it is large or small, though I think
    it would be better small. Of course I will share the expense
    willingly if it is heavy. \citep[15]{Russell2014}
\end{quote}

Perhaps we owe to the editorial largess of Morley the modern symbol
for disjunction, for, as Russell intimates, `$\lor$'~does not occur in
any of Whitehead's publications prior to Russell's 1906 piece.
Consequently, one reckons that Morley needed new type (of unknown
cost) for~`$\lor$'. This amusing episode in the history of~`$\lor$'
teaches us that editors should not scrimp on notations, and further,
that an erstwhile boon of co-authorship was sharing the cost of new
type.\footnote{For an in-depth discussion of the difficulties of
typesetting Frege's much more complex \emph{Begriffschrift} and the
importance of available type there, see
\citet{GreenRossbergEbert2015}.}

The crucial point, however, is that Russell expresses to Morley the
importance of not printing the disjunction sign as a round symbol. He
does not care about its size, but definitely wants it to be pointed. A
clear rationale for this is that Russell wants readers to readily
distinguish disjunction from class union, a distinction that Russell
had been concerned to make going back to \citet[Chapter
II]{Russell1903}. If the exact same sign was used for both, this could
confuse the reader concerning the logical distinction between
disjunction and union. This is why the size, in contrast, does not
matter, so long as the symbol is pointed: it has to be clearly
distinguishable on sight from its class theoretic analogue.

The 1906 paper accounts for the first published use of~`$\lor$' for
disjunction, but it is worth mentioning that the first unpublished
(and so the very first known) systematic uses of~`$\lor$' for
disjunction (and only for disjunction) occur somewhat earlier in
Whitehead and Russell's writings. Initially, Russell's practice was the reverse of
the modern one: he used the symbol~`$\lor$' for class union
and~`$\pmccup$' for ``propositional addition'' (disjunction). In his
spring 1902 ``On Likeness,'' \citep[440]{Russell1902} which was only
published after Russell's death, we find:\footnote{The definition may
be read as, ``$R$ and $R'$ are like relations (order isomorphic)''
means that there exists a relation in the intersection of one-one
relations and $S$s such that $S$'s domain~$\sigma$ is the union of the
domains of $R$ and $R'$, and is such that $R'$ is equal to the
relational product of $\pmcrel{S}$, $R$, and~$S$.''}
\[ \pmast1\pmcdot1\; (R)L(R') \pmdot = \pmdot R, R' \varepsilon
\text{Rel} \pmand \pmSome 1 \to 1 \pmand S \varepsilon (\sigma = \rho
\pmor \pmcrel{\rho} \pmand R' = \pmcrel{S} R S) \pmdf \] This
definition mashes notation from Peano's \emph{Formulaire} with some
symbols Russell independently devised in his 1901 ``The Logic of
Relations'' where he insisted on distinguishing relational products
and relational intersections \citep[315]{Russell1901a}. The key point
is that in his 1902 ``On Likeness'' manuscript, Russell is using the
symbol~`$\lor$' for the first time, although not for disjunction. Note
that he nonetheless distinguishes propositional disjunction from its
analogues for classes and relations. Again, this distinction could be
maligned or, as Peano would have it, exploited, but Russell is keen to
symbolically separate the two (formally different) notions.

Russell switches the symbols for disjunction and union in his early 1903 manuscript ``Classes,'' some of which is now lost. Here we find Russell's first known use of `$\lor$' for disjunction (and only for disjunction) \citep[9]{Russell1903b}:
\[
    \pmast12\pmcdot58\,\, \text{Quad}(\varphi) \pmdot\,=\,\pmdott (\pmSome f)\pmdot\{(x)\pmdott\varphi x \pmdot \pmiff \pmdot Fx\} \pmdot \pmor \pmdot (\pmSome f)\pmdot\{(x)\pmdott\varphi x \pmdot\pmiff\pmdot F'x\} \pmdf
\]
\begin{figure}
    \includegraphics[width=\textwidth]{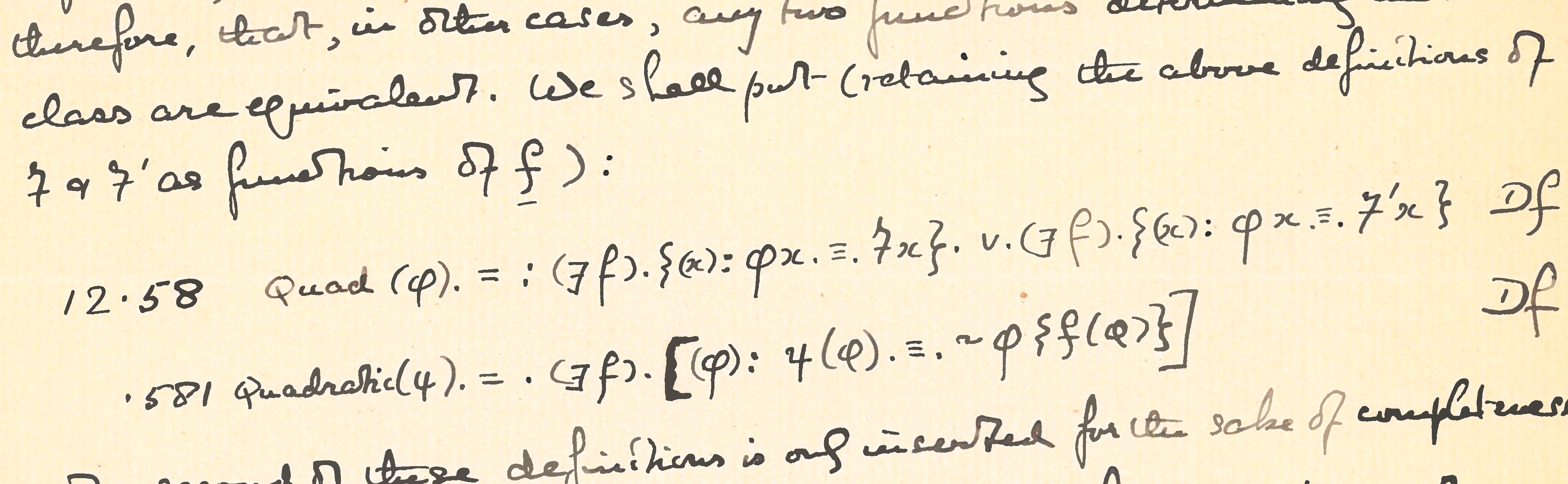}
    \caption{Use of $\lor$ for disjunction in \citet{Russell1903b}}
\end{figure}

Comments on a subsequent definition show that Russell here uses `$\lor$' \emph{only} for disjunction: 
\begin{quote}
    \[ \pmast14\pmcdot56\; a\pmccup b = x \mathrel{\varepsilon} (x
    \mathrel{\varepsilon} a \pmdot \pmor \pmdot x
    \mathrel{\varepsilon} b). \] If $a$ and $b$ are classes, $a\pmccup
    b$ is their logical sum, \emph{i.e.} the class formed by adding
    the two. If $x$ is either a member of $a$ or a member of $b$, then
    $x$ is a member of $a\pmccup b$, but not otherwise. The symbol
    `$a\pmccup b$' is read `$a$ or $b$'\dots\@ The function $a \pmccup
    b$ is the analogue for classes of $p \pmor q$ for propositions.
    \citep[19]{Russell1903b}
\end{quote}

Thus, `$\lor$' for disjunction was born in early 1903 (although
Russell used it for class union in 1902), and `$\lor$' for disjunction
made its first public appearance in April
1906.\footnote{Unfortunately, we cannot provide a photograph of the
manuscript page where Russell first defines `$p \pmor q$' because this
page is lost. We know it existed because \citet[20]{Russell1903b}
says, ``The definition of $p \pmor q$ in $\pmast4\pmcdot1$ cannot be
extended to a class of propositions\dots'' But at least what remains
of the manuscript allows us to give a literal picture of the earliest
known use of~`$\lor$' for disjunction in Russell's corpus.} This
answers the `when' part of question~(2). The `why' part, as we saw,
stemmed from Russell's desire to symbolically distinguish class union
and propositional disjunction.

There are two complications about this period between Russell's spring
1902 ``On Likeness'' (wherein the symbol `$\lor$' is first introduced)
and the 1903 ``Classes'' manuscript. First, as was pointed out to us
by Alasdair Urquhart, in a manuscript contemporaneous with the early
1903 ``On Classes,'' Russell uses `$\lor$' for disjunction (and only
for disjunction) as well. This manuscript consists of notes Russell
made on Frege's \emph{Grundgesetze} Vol. II and its attempted solution
to the contradiction---hence their title, ``Frege on the
contradiction.'' In one of the last definitions in this manuscript,
Russell defines class union using propositional sum, symbolized by
`$\lor$':
\[
     \pmcdot56 \quad a \pmccup b = x \mathrel{\rotatebox[origin=c]{180}{$\varepsilon$}} (x \mathrel{\varepsilon} a \pmdot \pmor \pmdot x \mathrel{\varepsilon} b) \pmdf  
\]
Urquhart's dating of this manuscript to is supported by a letter that 
Russell wrote to Frege, dated 20 February 1903, thanking him for sending 
Volume II of \textit{Grundgesetze} and asking about its solution to the 
contradiction. Given the content of these notes, it is likely they were 
written at nearly the same time as the ``On Classes'' manuscript. So we 
strictly-speaking have two candidates for the first use of `$\lor$' for 
disjunction---and imperfect information as to which one came first!
\begin{figure}
	\includegraphics[width=\textwidth]{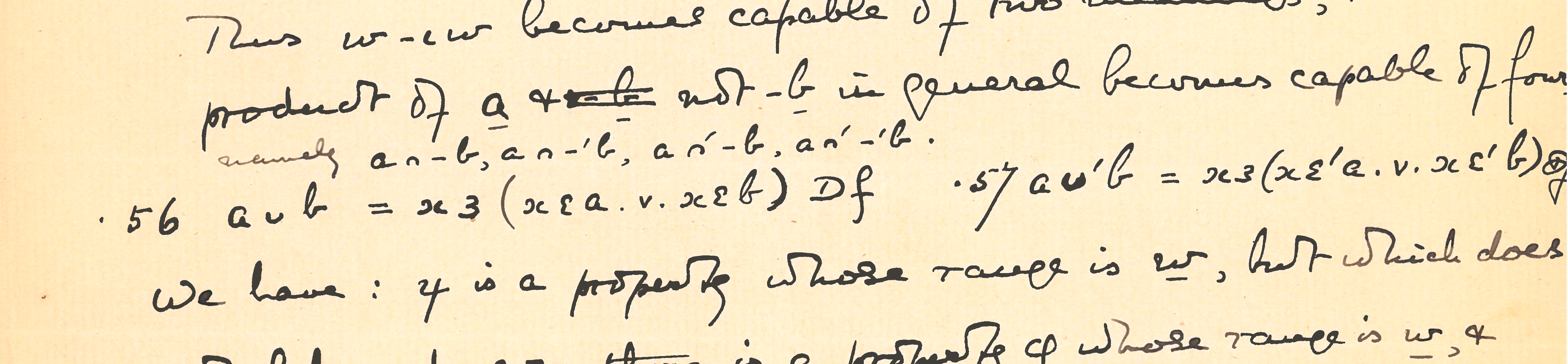}
	\caption{Use of $\lor$ for disjunction in \citet{Russell1903c}}
\end{figure}
The other complicating factor is that Whitehead also uses `$\lor$' for 
disjunction (and only for disjunction) in a manuscript that is also from 
early 1903, ``The Logic of Propositional Functions'' \citep{Whitehead1903}:
\begin{figure}
	\includegraphics[width=\textwidth]{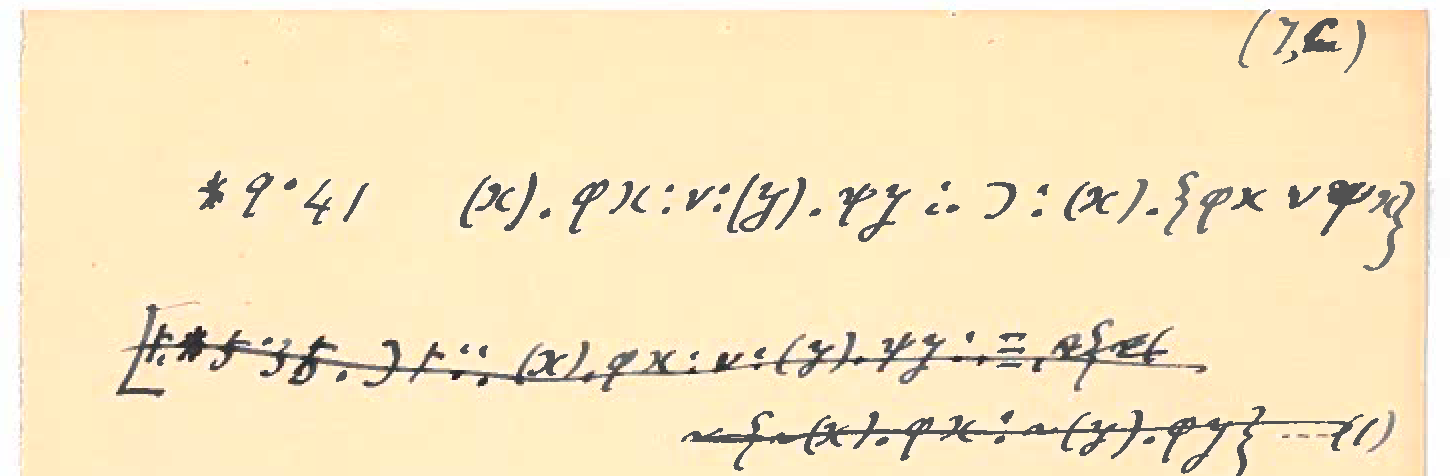}
	\caption{Use of $\lor$ for disjunction in \citet{Whitehead1903}}
\end{figure}
This is a familiar quantification theorem:
\[
     \pmast9\pmcdot41\quad (x)\pmdot \varphi x \pmdott \lor \pmdott (y)\pmdot \psi y \pmdot\hspace{.1em}\pmdott \pmimp \pmdott (x)\pmdot \{\varphi x \lor \psi x\}
\]
(In modern notation: $\forall x\,\varphi(x) \lor \forall y\,\psi(y)
\pmimp \forall x(\varphi(x) \lor \psi(x))$). It is not known for
certain when this manuscript was written, though the internal evidence
suggests it was either late 1902 or early 1903.\footnote{ See
Urquhart's headnote to \citet{Russell1903b}, p.~3 of
\citet{Russell1994b}.} Either way, Whitehead's manuscript dates from
after ``On Likeness'' wherein `$\lor$' is first used by Russell
(though not for disjunction). It is difficult beyond that to assign
priority for using `$\lor$' for disjunction. The earlier parts of
\citet{Russell1903b} were manuscripts on propositional logic and
likely included a definition of~`$\lor$' in terms of `$\pmimp$' and~$\pmnot$'. But it is not known whether Whitehead's or Russell's
manuscript came first. For this reason, we here give joint credit to
Whitehead and Russell.\footnote{One might argue that Russell is more
likely responsible because Russell in August 1902 drafted material
corresponding to \textit{Principia}'s propositional logic
($\pmast1$--$\pmast5$) \citep[691]{Moore1993}. Whitehead read this
material on propositional logic (and sharply criticized it). It would
not be surprising if Whitehead got the use of `$\lor$' for disjunction
from Russell's 1902 draft. Additionally, \citet[137-138]{Russell1948}
wrote, ``In the early part of \textit{Principia}, Whitehead
contributed the treatment of apparent variables and the notation $(x)
\pmdot \varphi x$. Chapters 10, 11, 13 of the \textit{Principia} are in
the main his work.'' It would be peculiar if Whitehead had contributed
notation to the propositional logic like `$\lor$' and Russell omitted
this point even while being sure to mention Whitehead's contribution
to the notation of \textit{Principia}'s quantifier theory.}

Russell's use of `$\lor$' for disjunction persists in Russell's
writings and through the 1910 publication of the co-authored
\emph{Principia}. The rest, as they say, is history: the symbol became
widely adopted following \emph{Principia}'s publication. In the next
section, we sketch how this state of affairs, which persists today,
came about.

As a coda to this part of the story, our suggestion that Russell 
chose `$\lor$' as a modification of~`$\pmccup$' is further supported 
by the surprisingly similar genealogy of `$\land$' and of dots for 
conjunction. Russell uses~`$\land$' for the first time in ``On 
Likeness'' for class intersection rather than for propositional 
product, as in the above definition $\pmast1\pmcdot1$.
But the story curiously diverges in the 1903 ``Classes'' manuscript:
there Russell introduces `$\land`m$' for a
class-operator---the conjunction of a class $m$ of propositions
($\pmast14\pmcdot8$)---and similarly introduces `$\pmor`m$'
for its dual class-operator---for the disjunction of a class $m$ of
propositions ($\pmast14\pmcdot81$). Despite this, even though Russell
also uses `$\lor$' for propositional disjunction, he uses Peano's
square dots for conjunction as well as scope by 1903. Russell was 
thus quite close to also giving us `$\land$' for conjunction, 
having used it briefly himself, but chose an alternative. 

Oddly, Russell seems not to have used `$\land$' for conjunction
anywhere else. Why did Russell not seize this opportunity? A partial
explanation is provided by the fact that the `$\land$' notation
``rarely occurs'' \citep[20]{Russell1903b}. So the `$\land$' symbol
was not much on Russell's mind even once he introduced it. Perhaps
more importantly, Russell had used dots for conjunction well before he
introduced either `$\lor$' or `$\land$.' For example, in his ``Logic
of Relations'' \citep{Russell1901e} Russell uses dots for conjunction
and for scope, reserving concatenation for relational products. In
contrast, Peano never used dots for conjunction: he only used dots for
scope \citep[118-119]{Schlimm2021}. \citet[\S1]{Peano1895} instead
uses concatenation for propositional conjunction and `$\pmccap$'~for
class intersection. So Russell perhaps inserted dots for conjunction
to distinguish conjunction from relational products. By the time he
thought to use `$\lor$' for anything, he had an established practice
of using dots for propositional conjunction. This squeezed out using
`$\land$' for that notion. This episode similarly shows that Russell
was keen to introduce new symbols to distinguish logical
notions---like propositional and relational ``conjunction''---and
further supports our hypothesis that this common feature of Russell's
notational design led to introduction of `$\lor$' rather than any love
for Latin or Leibniz.

\section{The adoption of `$\boldsymbol\vee$'}

The influence of \textit{Principia} cemented the use of the
Peano-Russell system of symbols among philosophers.  This is true
first of all of the logicians directly influenced by \textit{Principia},
including \citet{Wittgenstein1921} and \citet{Ramsey1926a}.
\citet{Carnap1929}, a simplified exposition of the type theory of
\textit{Principia}, uses the same notation, as does \citet{Quine1934}.
Early philosophy textbooks on logic that adopt a symbolic approach
(such as \cite{Lewis1918,NagelCohen1934,Langer1937}) also use the
symbols of \emph{Principia}. A notable exception, however, are the
Polish logicians of the Lw\'ow-Warsaw school, who used Polish notation
(usually $Apq$ for the alternation of $p$ and~$q$).  

Mathematicians working in the algebra of logic, such as
\citet{Lowenheim1915}, \citet{Skolem1920}, and \citet{Huntington1933},
continued to use the Boole-Peirce-Schr\"oder notation through the
1930s.\footnote{The use of and discussion of algebraic logic was not
confined to mathematicians. For instance, \citet{Mally1912} follows
Schr\"oder, but uses different symbols (and distinguishes between
propositional and class operations). \citet{Lewis1918} discusses the
``Boole-Schr\"oder algebra'' in the first chapter, as well as the
developments of \emph{Principia}. \citet{Jorgensen1931} faithfully and
exhausively describes the algebra of logic alongside the work of
Frege, Peano, and Whitehead and Russell, always using the same
symbolism as in the original works discussed. \citet{Langer1937} gives
equal space to the class calculus of Schr\"oder and \emph{Principia}:
``[T]his whole book has been planned to show very precisely that
`algebra of logic' and `logistic' are all of a piece\dots'' (p.~18).}
Beginning in the late 1920s, mathematicians more and more
converted to the notation made popular by
\emph{Grundz\"uge der theoretischen Logik}
\citep{HilbertAckermann1928}: overline for negation, `$\rightarrow$' for
the material conditional, `\&' for conjunction---and `$\lor$' for
disjunction. 

Hilbert and many of his students had mainly relied on
\citet{Schroeder1890} as their source of logic, until 1918, when
Hilbert tasked Heinrich Behmann with investigating the new approach to
the foundations of mathematics of \emph{Principia} (see
\cite{Behmann1918,Mancosu1999}). In his 1917/18 and 1920 lectures on
logic \citeyearpar{Hilbert1918,Hilbert1920a}, Hilbert used `$+$' and
`$\times$' (or just juxtaposition) for conjunction and disjunction (with
overline for negation)---the opposite of the way `$+$' and~`$\cdot$' are
used in the algebraic tradition. \citet{Behmann1922,Behmann1927} was
the only one to use this notation in print.\footnote{In fact, Behmann
expressed everything using only disjunction (symbolized as
juxtaposition) and negation (overlining).} However, he did this for
the specific purpose of logical calculation, where compactness and
ease of manipulation is important. He generally endorsed the symbolism
of \emph{Principia} on grounds of uniformity, writing at the beginning
of \citet{Behmann1922}, 
\begin{quote}
    One of the most important requirements for the following
    investigation is that we provide a suitable symbolism. It may seem
    to be a given that we should replace the old symbolism of
    Schr\"oder by the newer one which is currently the most highly
    regarded and seems to have the best chance of becoming widely
    adopted, namely the one of \emph{Principia Mathematica} by
    Whitehead and Russell.\footnote{``\foreignlanguage{german}{Eines
    der wichtigsten Erfordernisse f"ur die folgende Untersuchung ist
    die Bereitstellung einer geeigneten Symbolik. Es m"ochte hier als
    das Gegebene erscheinen, die alte Schr"odersche Symbolik durch
    diejenige neuere zu ersetzen, die gegenw"artig am meisten
    gesch"atzt wird und wohl am ehesten Aussicht hat, allgemein
    durchzudringen, n"amlich die der Principia Mathematica von
    Whitehead und Russell}.'' Our translation.}
\end{quote}
The treatment of logic in both Hilbert's lectures and in Behmann's
work is already axiomatic (as in \emph{Principia}), and not algebraic.

The idiosyncratic notation of 1920 was replaced, in the lecture notes
to the 1921/22 course ``Foundations of mathematics''
\citeyearpar{Hilbert1921}, by the Hilbert-Ackermann notation, which
used the symbols `$\lor$', `\&', `$\to$', and overlining for negation. A
handwritten note by Bernays in the copy of \citet{Hilbert1920a} filed
in the Göttingen mathematics library, headed ``Suggestions for
nomenclature [\emph{\foreignlanguage{german}{Vorschläge zur
Bezeichnung}}],'' reads in part:
\begin{quote}
$X \lor Y$: $X$ or $Y$ (abbreviation) (``$\lor$'' is Russell's symbol,
v initial letter of ``vel'')
\footnote{``\foreignlanguage{german}{$X \lor Y$: $X$ oder $Y$
(Abkürzung) ("`$\lor$"' Zeichen von Russell, v Anfangsbuchstabe von
"`vel"')}'' \citep[334]{Hilbert2013}. On p.~2, where conjunction and
disjunction are introduced as $X + Y$ and $X \times Y$, respectively,
marginal annotations record the symbols later preferred: $X
\mathbin{\&} Y$, $X \lor Y$ \citep[300]{Hilbert2013}.}
\end{quote}
It thus seems likely that the motivation to adopt `$\lor$' for
disjunction \emph{in the Hilbert school} was two-fold: first, to
follow Russell's established usage, and second, because the symbol is
reminiscent of \emph{vel}.

The first use of this system of logical notation in print is in
\citet{Hilbert1923}, and it was adopted in all subsequent publications
on logic in Hilbert's school.\footnote{\citet{Gentzen1935a} is an
exception: he used Russell's `$\pmimp$' for the conditional and
`$\to$' for the sequent arrow. For more on the development of
propositional logic in the Hilbert school, see \citet{Zach1999}.}
\citet{HilbertAckermann1928} proved extremely influential. With it,
Hilbert's way of writing formulas became widely adopted, especially
after \citet{Godel1930,Godel1931}---even Skolem used it exclusively
from the late 1930s onward (e.g., see \citealt{Skolem1935a}).

Tarski was catholic in his use of notations: he used whatever seemed
most appropriate. In his dissertation on the logic of
\emph{Principia}, he applied the notation of \emph{Principia}
\citep{Tajtelbaum1923}. Writing on propositional calculi with \L
ukasiewicz, he used Polish notation \citep{LukasiewiczTarski1930}.
Writing for an audience of set theorists, he used algebraic notation
\citep{KuratowskiTarski1931}. Writing for an audience of philosophers
(in \emph{Erkenntnis}, the journal of the Vienna Circle), he used the
notation of \emph{Principia} \citep{Tarski1934}. And when he wrote for
his mathematical colleagues in Vienna (Hahn, G\"odel), his logical
language used the Hilbert-Ackermann symbols
\citep{Tarski1933}.\footnote{Notably, \citet{Tarski1933} also used
$\land$ for conjunction; see below.} 

\citet{Bernays1926} mentions that `$\lor$' is to have the inclusive
meaning of the Latin ``\emph{vel}.'' \citet{BernaysSchonfinkel1928}
and \citet[p. 4]{HilbertAckermann1928} likewise mention this, and also
note that~`$\lor$'---typeset in the latter as a lowercase~`v'---is not
to be confused with exclusive or ``in the sense of the latin
`aut--aut'.'' The identification of \textit{vel} with inclusive and
\textit{aut} with exclusive or is already made in
\citet[p.~226]{Schroeder1890}, who however used `$+$' and
never~`$\lor$'.\footnote{\citet{Jennings1994} argues that it is a myth
that \textit{aut} has the sense of exclusive disjunction.
Incidentally, \citet{Schroeder1890} discusses Leibniz extensively, but
only cites the Erdmann edition \citep{Leibniz1840} as well as
\citet{Kvet1857}, neither of which mention `v' as a symbol for
disjunction or union.}

Quine's \emph{Mathematical logic} of 1940 may be the source of the
connection between `$\lor$' and the initial letter of \textit{vel} in
English textbooks.\footnote{Other early textbooks in English that use
symbolic logic \citep{Lewis1918,NagelCohen1934,Langer1937,Tarski1941a}
do not make this connection even though they use the \emph{Principia}
notation.} He first discusses the exclusive and inclusive use of `or'
in English:
\begin{quote}
We must decide whether `or' is to be construed in an \emph{exclusive}
sense, corresponding to the Latin `aut', or in an \emph{inclusive}
sense, corresponding to the Latin `vel'. \dots When `or' is used in
the inclusive sense, on the other hand, the compound is regarded as
true if at least one of the components is true; joint truth of the
components verifies the compound. 
\dots\@
In mathematical logic the ambiguity of ordinary usage is resolved by
adopting a special symbol `$\lor$', suggestive of `vel', to take the
place of `or' in the inclusive sense. Alternation is identified with
this usage \dots\@ The exclusive use of `or' is not frequent enough in
technical developments to warrant a special name and symbol.
\citep[12--13]{Quine1940}
\end{quote}
From \emph{Mathematical Logic}, the idea that \emph{aut} and
\emph{vel} in Latin represent the exclusive and inclusive sense of
`or', and the use of the symbol~`$\lor$' for the latter, percolates
into other introductory logic texts, first Quine's own \emph{Methods
of logic}:
\begin{quote}
Latin has distinct words for the two senses of `or': \emph{vel} for
the nonexclusive and \emph{aut} for the exclusive. In modern logic it
is customary to write `$\lor$', reminiscent of \emph{vel}, for `or' in
the nonexclusive sense: `$p \lor q$'. It is this mode of compounding
statements, and only this, that is called \emph{alternation}.
\citep[5]{Quine1950a}
\end{quote}
Quine adds a section on ``Alternation and duality'' to the revised
edition of \emph{Elementary logic}:
\begin{quote}
Special symbols are commonly added for `or' and `if'. For `$p$ or $q$'
the notation is `$p \lor q$'; here `$\lor$' stands for the Latin
\emph{vel}, which means `or' in the inclusive sense.
\citep[52]{Quine1965}.\footnote{The first edition \citep{Quine1941}
symbolized everything with only negation `$\sim$' and conjunction
`$\cdot$' and introduced alternation as the denial of `neither\dots
nor'.}
\end{quote}
As we saw above, though, Quine's alleged connection between `$\lor$'
and \emph{vel} is likely spurious: given the textual evidence that we
have, the only way this conjecture would hold good would be if Russell
saw Peano's mention of the ``$\lor$~is for \textit{vel}'' story about
Leibniz and, for that reason, Russell decided to use that notation
himself. We argued above that this was a very remote possibility for 
which there is no surviving evidence, and that there is a highly 
plausible alternative explanation that fits with a common feature 
of Whitehead and Russell's notational design practices. 

The work of Whitehead and Russell, and of Hilbert and his students,
was defining for the further development of logic in the 20th century.
This gives us an answer to the further question, ``why do \emph{we}
(still) use~`'$\lor$'?'' The answer is that we have inherited our logical
notations from two influential research projects, namely that of
\emph{Principia} and the project of formalizing first-order logic and
investigating its properties (completeness, undecidability, etc.)
begun by Hilbert.
Both of these used `$\lor$' for disjunction.  Textbook presentations
cemented this use further, and Quine's texts were likely mainly responsible
for the spread of the myth that `$\lor$' abbreviates
\emph{vel}.\footnote{This is not to say that alternative notations
were not used after, say, 1940. Polish notation, in particular,
remained popular through much of the 20th century.}

\section{Other symbols}

The symbol `$\lor$' is special in that it now is almost universally
used for disjunction. But it has been used for other purposes as well;
and other symbols also have interesting stories.

Ladd-Franklin \citep[25--26]{Ladd1883} used `$\lor$' and `$\laddexcl$'
in her dissertation, but only for the copula, that is, for
predications among classes and not for propositional connectives. So,
Ladd says that `$A \lor B$' means ``$A$ is in part~$B$,'' that is,
``some $A$ is~$B$,'' and `$A \laddexcl B$' means ``$A$ is-not~$B$,'' or
``No $A$ is~$B$.'' The uses of `$\lor$' and `$\laddexcl$' thus
correspond to, as Ladd puts it, (Boolean class) ``inclusions'' and
``exclusions'' \citep[26]{Ladd1883}. She deals, rather, with an
analogue of propositional disjunction from a Boolean perspective,
where the symbol occurs between categorical terms or concepts rather
than between truth-apt
formulas.\footnote{See \citet{Uckelman2021} for a discussion of
Ladd-Franklin's dissertation, and \citet{Wege2020} for a discussion of
Ladd-Franklin's hypothesis that commutative operators should be
represented by symmetric symbols.}

\citet{Tarski1933} was one of the first to combine `$\lor$'
with~`$\land$' for conjunction. \citet{Heyting1930} predates this use,
however, and is perhaps the first to use `$\land$' in print as the
symbol for conjunction. The first use of~`$\land$' for conjunction in
a textbook is probably \citet{Tarski1941a} (the original Polish and
German editions did not use logical symbols).

\citet{Sheffer1913} used `$\land$' as the symbol for ``neither \dots
nor'' (\textsc{nor}), i.e., the logical connective
\emph{corresponding} to the operator in Boolean algebra which
\emph{he} symbolized~`$\boldsymbol\mid$'. \citet{Nicod1917} used `$\boldsymbol\mid$' instead
for ``not both'' (\textsc{nand}), and gave a reduction of the logical
constants of \emph{Principia} to it (the reason being that the
resulting definitions of other connectives turn out to be much simpler
than with Sheffer's \textsc{nor}). This has stuck; ``Sheffer stroke''
now universally means \textsc{nand}. We should point out here,
however, that the first published discussion of \textsc{nand} and
\textsc{nor}, including proofs that the usual primitives can be
expressed by either, and an axiomatization of Boolean algebra in terms
of them, is neither \citet{Sheffer1913} nor \citet{Nicod1917}, but
\citet{Stamm1911}.\footnote{Stamm used $\invlazys$ for \textsc{nand}
and $\ast$ for \textsc{nor}. His paper seems to have excited little
interest at the time, or since. L\"owenheim belittled the results as
``not especially useful'' in his review in the \emph{Jahrbuch f\"ur
die Fortschritte der Mathematik}~42.0079.02.}

The ``neither \dots nor'' connective (\textsc{nor}) is now commonly
called the Peirce arrow, and written~`$\boldsymbol\downarrow$'. However, Peirce
himself did not use `$\boldsymbol\downarrow$' for \textsc{nor}, rather, it was
introduced by \citet{Quine1940}:
\begin{quote}
The fact that $\ulcorner(\phi \peircearrow \psi)\urcorner$ denies
$\ulcorner(\phi \lor \psi)\urcorner$ is reflected in ordinary
language, indeed, by the cancellatory `n' which turns `either---or'
into `neither---nor'. The vertical mark in `$\peircearrow$' may be
thought of as having the same cancellatory effect upon the
sign~`$\lor$'; it is like the vertical mark in the inequality sign
`\ooalign{\hfil$\boldsymbol=$\hfil\cr\hfil$\boldsymbol|$\hfil}' of arithmetic. (p.~46)
\end{quote}
Although Quine himself does not call his `$\peircearrow$' the `Peirce
arrow', he does credit Peirce for his work on it in fn.~1 on p.~49: 
\begin{quote}
The definability of denial, conjunction, and alternation in terms of
joint denial [\textsc{nor}, the Peirce arrow] was first pointed out by
Sheffer in \citeyear{Sheffer1913}; and similarly for alternative
denial [\textsc{nand}, the Sheffer stroke]. The adequacy of joint
denial was known to Peirce in 1880, and both facts were known to him
in 1902; but his notes on the subject remained unpublished until
\citeyear{Peirce1933} (4.12, 4.264).
\end{quote}
As Quine remarks, Peirce's writings on \textsc{nor} had only recently
become available. In the manuscripts mentioned by Quine, Peirce first
simply used juxtaposition: ``Two propositions written in a pair are
considered both denied. Thus, $AB$ means that the propositions $A$ and
$B$ are both false'' \citep[13]{Peirce1880}. In \citet{Peirce1902}, he
introduced the symbol~`$\ampheck$', which he called ``ampheck''
(p.~215).

\citet{Quine1940} also gave an origin story for~`$\sim$':\footnote{In
\citet{Quine1950a}, he prefers `$-$' for negation, and~`$\bar p$' for
negated atoms. Incidentally, \citet{Quine1950a} also seems to be the
origin of `$\top$' and~`$\bot$' as constants for truth and falsity:
``A convenient graphic method of imposing interpretations, of the
second of the above varieties, is simply to supplant the letters in a
schema by the mark `$\top$' for truths and `$\bot$' for falsehoods''
(p.~23).}
\begin{quote}
In mathematical logic the denial of a statement is formed by prefixing
the tilde `$\sim$', which is a modified `n' and is conveniently read
`not'. (p.~14)
\end{quote}
The history recounted in this paper shows that we should be wary of
just-so stories about the origins of notations like Quine tells here.
If we had to make a conjecture, we would agree with
\citet[lxi]{Moore2014} that Russell began using a tilde for negation
for reasons that conform to his well-attested to practice of
distinguishing different notions in his symbolism whenever the formal
laws differed, so that `$\pmnot$' was introduced to distinguish
propositional negation from class (or arithmetic) subtraction. 

The two standard symbols for the material conditional (implication)
are `$\pmimp$' and `$\to$'. The latter, as noted above, was
introduced by \citet{Hilbert1918,Hilbert1920a} and become widely
adopted after it was used in the influential \emph{Grundzüge} textbook
\citep{HilbertAckermann1928}. The origin of~`$\pmimp$' is a bit more
complicated: Peano introduced~`$\conseq$' as a symbol for both the
conditional and for class containment, as the rotated version of
`\textsc{c}', abbreviating ``\emph{est consequentia}'' and
``\emph{continet}'', Latin for ``is a consequence of'' and
``contains'':
\begin{quote}
    [The sign \textsc{c} means \emph{is a consequence of;} thus $b
    \mathbin{\textsc{c}} a$ is read \emph{$b$~is a consequence of the
    proposition~$a$}. But we never use this sign.]

    The sign $\conseq$ means \emph{one deduces}; thus $a \conseq  b$
    means the same as $b \mathbin{\textsc{c}}
    a$.\footnote{``\foreignlanguage{latin}{[Signum \textsc{c} significat
    \emph{est consequentia;} ita $b \mathbin{\textsc{c}} a$ legitur
    \emph{$b$~est consequentia propositionis~$a$}. Sed hoc signo
    nunquam utimur.]\par Signum $\conseq$ significat \emph{deducitur;}
    ita $a \conseq b$ idem significat quod $b \mathbin{\textsc{c}}
    a$}'' (\cite[viii]{Peano1889}, translation from
    \cite[105]{Kennedy1973}).}

    The sign $\conseq$ means \emph{is contained.} Thus $a \conseq  b$
    means \emph{the class~$a$ is contained in the class~$b$}. \dots
    The signs $\Lambda$ and~$\conseq$ have meanings here which are
    slightly different from the preceding, but no ambiguity will
    arise, for if propositions are being considered, the signs are
    read \emph{absurd} and \emph{one deduces}, but if classes are
    being considered, they are read \emph{empty} and \emph{is
    contained}.\footnote{``\foreignlanguage{latin}{Signum $\conseq$
    significat \emph{continetur}. Ita $a \conseq  b$ significat
    \emph{classis $a$ continetur in classis~$b$}. \dots Hic signa
    $\Lambda$ et $\conseq$ significationem habent quae paullo a
    praecedenti differt; sed nulla orietur ambiguitas. Nam si de
    propositionibus agatur, haec signa legantur \emph{absurdum} et
    \emph{deducitur;} si vero de classibus, \emph{nihil} et
    \emph{continetur}}'' (\cite[xi]{Peano1889}, translation from
    \cite[108]{Kennedy1973}).}
\end{quote}
Beginning with \citet{Peano1898}, the symbol~`$\conseq$' is replaced
with~`$\pmimp$' but with the same ambiguous use.\footnote{As with
Peano's~`$\pmccup$' and Leibniz's~`v', one may wonder if Peano came up
with~`$\conseq$' on his own or if he adopted someone else's use.  In
the \emph{Formulaire}, Peano was very generous with giving credit and
listing the symbols used by others.  However, of `$\varepsilon$'
and~`$\conseq$' he wrote, in \citet[29]{Peano1897}, that they are
symbols that he introduced, just before he surveys symbols others used
for the same purpose. Only sometime around 1899 did he became aware of
prior uses of~`$\conseq$'. He mentions a use of~`$\conseq$' in the
sense of ``therefore'' between theorems by \citet[36]{Abel1881} in
\citet[8]{Peano1899}. In \citet{Peano1903}, he also mentions
\citet{Gergonne1816}, who used `$C$'
and~`\raisebox{\depth}{\rotatebox{180}{$C$}}' for ``contains'' and
``is contained in.'' Gergonne does not use them as propositional or
class operators, however, but markers of different kinds of
subject-predicate propositions. \citet[129]{Schroeder1890}
independently introduced signs similar to `$\boldsymbol\subset$'
and~`$\pmimp$' as signs for proper concept containment in the same
directions we now use them (e.g., ``gold $\boldsymbol\subset$
metal''), and uses~`$\subsumed$' for contains-or-equals. His
comparison with the `$<$' symbol suggests he also was following the
design principle that analogous notions should be symbolized with
similar symbols.} Russell of course derived his~`$\pmimp$' from Peano
and used it as early as \citep{Russell1901e}. This was a French
translation of a 1900 manuscript \citep{Russell1900b}, published in
Peano's \emph{Revue de mathématiques}. Therein the symbol was typeset
interchangeably as `$\pmimp$' and as~`$\conseq$'.\footnote{Russell was
displeased with the typesetting of this article, as he complained to
Couturat in a letter of 27 July 1901; see the headnote to
\citet{Russell1901e}. Curiously, the mishmash of `$\pmimp$'
and~`$\conseq$' is still present in the collection
\citet[1--38]{Russell1956}, despite Russell correcting and approving
that English translation of \citet{Russell1901e}. Note that the
reprints in \citep{Russell1993} correct these misprints, and render
implication uniformly as~`$\pmimp$'.}  Later, including in
\textit{Principia}, Russell had it typeset as~`$\pmimp$'.

\section{Conclusion}

Logicians and historians of logic have suggested various origin
stories for~`$\lor$' as the symbol for disjunction. As we saw above,
Peano suggests that Leibniz chose `$\lor$' for the Latin \emph{vel}.
Kneale and Kneale suggest that Peano first chose it. Cajori attributes
the first use to \emph{Principia}. According to Quine, `$\lor$'~stands
for, or at least suggests, the Latin \emph{vel}. The evidence
discussed in this paper indicates that `$\lor$'~for disjunction is not
a straightforward Latin abbreviation, and that logicians and
historians of logic have often attributed the first use of `$\lor$'
for disjunction to the wrong author or to the wrong text.

The appealing story often told that our `$\lor$' for disjunction comes
from the Latin \emph{vel} is not well-supported by the textual record.
It was Whitehead and Russell who first systematically used~`$\lor$' 
for disjunction, and it is more likely given Whitehead and Russell's 
notational habits (and Russell's disgust for Latin) that this choice of 
notation was to stress the analogy of propositional addition with class union.

This may all seem like nothing more than logical \emph{Trivial
Pursuit}. However, our discussion also highlights features of the
trajectory of the development of logical notations that is not
trivial.  We have seen that the considerations for or against the
choice of a symbol used for `or' were often deeper than just whether
or not it serves as a good mnemonic device (a symbol reminiscent of
\emph{vel}). The mathematical tradition of the algebra of logic
stressed the similarity of class union and propositional (or
truth-value) disjunction with algebraic operations, and thus often
chose the same symbols used in algebra (`$+$' or~`$\times$'). Other
writers were concerned with separating arithmetical or algebraic
notations (or defining them) from more general logical notations.
This, e.g., led Peano to use a \emph{different} symbol for (logical)
union and disjunction than the algebraic~`$+$', and it led Whitehead
and Russell to use distinct symbols for class union and propositional
disjunction. These concerns for mathematical uniformity and
philosophically relevant distinctions, of course, extend to symbols
and notions other than~``or.''

The story connects with other important concerns that influenced the
development of a system of notation which, since Cajori penned the
passage we opened the paper with, has become \emph{almost}
standardized. These include the motivations for choosing the
particular set of primitives that have become commonplace, e.g., why
we take inclusive and not exclusive disjunction as a primitive. It
also connects with the question of why and when different approaches
to logic (and notation) became more widespread or lost appeal: the
choice of notations of the pioneering writers of the early to mid-20th
century provides insight into their influences and methodological
commitments. We have only touched upon these deeper topics, but they
suggest fruitful directions for further research.

If nothing else, we now have an answer to the question:
``Why~$\lor$?'' We owe it to Whitehead and Russell, who in spring 1902
wanted a symbol different from Peano's (and Grassmann's) symbol
`$\pmccup$' to distinguish union and propositional disjunction in his
symbolism. After, and in part because Hilbert endorsed Whitehead and
Russell's choice to use `$\lor$' for disjunction in
\textit{Principia}, the symbol became ubiquitous.

\subsection*{Acknowledgments} 

We thank the hosts Alexander Klein and Sandra Lapointe as well as the
participants of the Fall 2020 Bertrand Russell Research Centre and
\emph{Journal of the History of Analytical Philosophy} graduate
seminar, who raised insightful questions and gave helpful comments.
The paper was also discussed with the History and Philosophy of
Mathematics Reading Group at McGill University, and we especially
thank Moritz Bodner, Julien Ouellette-Michaud, David Waszek, and Dirk
Schlimm for their questions and suggestions. We are also grateful to
Katalin Bimb\'o, Edward Buckney, Bernard Linsky, Paolo Mancosu, Marcus
Rossberg, Sara~L. Uckelman, Rafa\l{} Urbaniak, and especially Massimo
Mugnai, Brent Odland, Alasdair Urquhart, Jack Zupko, and the reviewer
for the~\emph{RSL}, all of whose suggestions and comments improved the
manuscript. We'd also like to acknowledge the Bertrand Russell
Archives at McMaster University and the Gottfried Wilhelm Leibniz
Bibliothek for providing scans of the Whitehead, Russell, and Leibniz
manuscripts. Elkind thanks the Killam Trusts for their support.

\nocite{Kennedy1973}

\printbibliography


\vspace*{10pt}
\noindent
WESTERN KENTUCKY UNIVERSITY\\
\hspace*{9pt}DEPARTMENT OF POLITICAL SCIENCE\\
\hspace*{18pt}1906 COLLEGE HEIGHTS BLVD \#31086\\
\hspace*{27pt}BOWLING GREEN, KY 42104, USA\\
{\it E-mail}: \href{mailto:landon.elkind@wku.edu}{\texttt{landon.elkind@wku.edu}}\\
{\it Webpage}: \url{https://landonelkind.com/}\\[8pt]
UNIVERSITY OF CALGARY\\
\hspace*{9pt}DEPARTMENT OF PHILOSOPHY\\
\hspace*{18pt}2500 UNIVERSITY DR NW\\
\hspace*{27pt}CALGARY, AB T2N 1N4, CANADA\\
{\it E-mail}: \href{mailto:rzach@ucalgary.ca}{\texttt{rzach@ucalgary.ca}}\\
{\it Webpage}: \url{https://richardzach.org/}

\end{document}